\documentclass[11pt]{article}
\usepackage{graphicx}
\usepackage{latexsym}
\usepackage{amsmath}
\usepackage{amssymb}
\usepackage{hyperref}
\usepackage{amsfonts}
\usepackage{algorithm}
\usepackage{algorithmic}
\usepackage{color,xcolor}
\usepackage{manfnt}
\usepackage{float}

\newcommand{\qed}{\nobreak \ifvmode \relax \else
      \ifdim\lastskip<1.5em \hskip-\lastskip
      \hskip1.5em plus0em minus0.5em \fi \nobreak
      \vrule height0.75em width0.5em depth0.25em\fi}

\def \Vh0{\stackrel{\circ}{V}_h} 
   
\def\Om{\Omega}   
\newcommand{\q}{\quad}

  \def\x{{\bf x}}

\def\bb{\begin{equation}} \def\ee{\end{equation}}
\def\beqn{\begin{eqnarray}}  \def\eqn{\end{eqnarray}}
\def\beqnx{\begin{eqnarray*}} \def\eqnx{\end{eqnarray*}}

\title{A parallel space-time domain decomposition method
for unsteady source inversion problems}
\author{Xiaomao Deng \footnote{
Laboratory for Engineering and Scientific Computing,
        Shenzhen Institutes of Advanced Technology, Chinese Academy of Sciences,
        Shenzhen, Guangdong 518055, {P. R. China}. The work of this author was partly supported by NSFC grant 91330111 and 11501545.
        (xm.deng@siat.ac.cn)}
\quad\quad Xiao-chuan Cai\footnote{Department of Computer Science, University of Colorado Boulder, Boulder, CO 80309, {USA}.  The work of this author was partly supported by NSF grant CCF-1216314.
(cai@colorado.edu)}
\quad\quad Jun Zou\footnote{Department of Mathematics, The Chinese University of Hong Kong, Shatin N.T., Hong Kong, {P. R. China}. The work of this author was substantially supported by
the Hong Kong RGC grants (Projects 405513 and 404611). (zou@math.cuhk.edu.hk)}}

\date{}
\begin{document}
\maketitle
\begin{abstract}
In this paper, we propose a parallel space-time domain decomposition method for solving an unsteady
source identification problem governed by the linear convection-diffusion equation.
Traditional approaches require to solve repeatedly a forward parabolic system, an adjoint system and a
system with respect to the unknowns.
The three systems have to be solved one after another.
These sequential steps are not desirable for large scale parallel computing.
A space-time restrictive additive Schwarz method is proposed for a fully implicit space-time coupled discretization scheme
to recover the time-dependent pollutant source intensity functions.
We show with numerical experiments that the scheme works well with noise in the observation data.
More importantly it is demonstrated that the parallel space-time Schwarz preconditioner is scalable on a supercomputer with over $10^3$ processors, thus promising for large scale applications.

\end{abstract}

\section{Introduction}
Pollutant source inversion problems have wide applications in, for example, the detection and monitoring of indoor and outdoor air pollution, underground water pollution, etc.
In the last several decades, physical, chemical and biological technologies have been developed  to identify different types of sources \cite{AB01,ZZFZG11,ZCP01}.
In this paper, assuming the pollutant concentration data is measured by distributed sensors, we
reconstruct the source intensities numerically using noise-contaminated data.
Like all inverse problems, such a reconstruction problem is ill-posed in the sense of Hadamard \cite{EHN98,SV07,W03}.
The lack of stability with respect to the measurement data is a major issue, which means that small noise in the data may lead to significant changes in the reconstructed source strength.
This problem has attracted much attention, and various methods have been developed, including both deterministic and statistical methods \cite{LZ07,SK97}.
Among the deterministic methods, quasi-explicit reconstruction formulas are available for one-dimensional source location recovery problems \cite{H09, H092};
and quasi-reversibility methods can be used to retrace the pollutant history as in \cite{SK95};
optimization based methods are also widely used \cite{AB03,KZ98,KZ00,SK94,WY10,XZ05}.
By reformulating an inverse problem into an output least-squares PDE-constrained optimization problem complemented with Tikhonov regularization, classical optimization methods such as
regression methods \cite{GER83}, linear and nonlinear programming methods \cite{GER83}, linear and nonlinear conjugate gradient methods \cite{AB05,WY10}, Newton type methods, etc. can be used to obtain the approximate solutions.
These methods can be categorized as sequential quadratic programming (SQP) methods.
Reduced space SQP methods decouple the system and iteratively update the state variable, the adjoint variable and the optimization variables by solving each subsystem in a sequential order.
In some sense this is a block Gauss-Seidel iteration with three large blocks.
Such methods require less memory due to the reduced subproblem size but the number of outer iterations for a specified accuracy grows quickly with the increase of the optimization variables,
thus they are not ideal for supercomputers with a large number of processors.
We introduce in this paper a full space approach that does not have the three large sequential steps as in the reduced space approaches.
Similar approaches have been applied to flow control problems in \cite{PBC06}.
The full space method solves the state variable, adjoint variable and the optimization variables simultaneously, thus avoids repeatedly solving the subsystems.
However the fully coupled system is several times larger in size and more ill-conditioned,
direct methods such as Gaussian elimination or LU factorization as well as the classical iterative methods such as the Jacobi method, the Gauss-Seidel method are not suitable.
To ease the difficulty of solving the large system, a preconditioned Krylov subspace technique is considered to reduce the condition number and the computing time significantly \cite{CLZ09, YPC12}.

The inverse problem of recovering the pollutant source intensity functions can be reformulated into a PDE-constrained optimization problem.
In this paper, we derive its continuous Karush-Kuhn-Tucker (KKT) system \cite{KT51}, including the state equation, the adjoint equation and other derivative equations with respect to each unknown source intensity.
Two main challenges of the problem lie in that firstly the adjoint equation needs the final state of the pollutant source distribution,
which implies that the state equation and the adjoint equation should be solved in a sequential order;
secondly the time marching of the unsteady problem is directional and sequential, thus difficult to break down into parallel steps.
For unsteady PDE-constrained optimization problems, a steady state optimization subproblem is solved at each time step \cite{YPC12}.
And in \cite{H08}, a block time-marching method is used to reduce the number of sequential steps and increase the degree of parallelism.
In this paper, we propose a fully coupled space-time domain decomposition method that couples the time with the space domain and decomposes the ``space-time'' domain into sub-domains,
then apply an additive Schwarz preconditioned Krylov subspace technique to solve the ``space-time'' problem.
Our algorithm is fully parallel in space and time, avoids the sequential time marching steps, and does not need to solve optimization subproblems.
As far as we know, no published work has achieved such a degree of parallelism for time-dependent inverse problems.

The rest of this paper is arranged as follows. The mathematical model and its corresponding optimization functional, and the derivation of the KKT system are formulated in Section \ref{sec:model}. 
The discretization of the KKT system is given in Section \ref{sec:finite element}.
The parallel algorithm for solving the KKT system is proposed in Section \ref{sec:scalable solvers}.
Some numerical experiments are shown in Section \ref{sec:numerical examples} and concluding remarks are given in Section \ref{sec:concluding remarks}.

\section{Model formulation}\label{sec:model}
We consider a flow domain $\Omega\in \mathbf{R}^2$ in which several point pollutant sources are present.
The distribution of the pollutant concentration is denoted by  $C(\mathbf{x},t)$ at location $\x$ and time $t$.
The transport process is modeled by the following convection-diffusion equation \cite{AB01,RR05}:
\begin{equation}
{\,\ } \hskip0truecm
\frac{\partial C}{\partial t}=\nabla\cdot(a(\mathbf{x})\nabla C)-\nabla\cdot (\mathbf{v}(\mathbf{x})C)+\sum_{i=1}^s \delta(\mathbf{x}-\mathbf{x}_i^*)f_i(t), ~0< t< T, ~\mathbf{x}\in \Omega,
\label{eq:cfeq}
\end{equation}
where $f_i(t)$ is the temporal intensity of the $i^{\mbox{\tiny th}}$ source at location $\mathbf{x}_i^*$, $i=1,\cdots,s$, $s$ is the number of sources, $a(\mathbf{x})$ and $\mathbf{v}(\mathbf{x})$ are the diffusive and convective coefficient.
$\delta(\cdot)$ is the Dirac delta distribution \cite{B03}.
The model is complemented by  the following boundary conditions
\begin{equation}\label{eq:BCs}
\centering{
C(\mathbf{x},t)=p(\mathbf{x},t),\quad \x\in \Gamma_1; \q
a(\mathbf{x})\frac{\partial C}{\partial \mathbf{n}}=q(\mathbf{x},t),\quad \x\in \Gamma_2
}
\end{equation}
and the initial condition
\begin{equation}\label{eq:IC}
C(\mathbf{x},0)=C_0(\mathbf{x}), \q \x\in \Om,
\end{equation}
where $\Gamma_1$ and $\Gamma_2$ cover the physical boundary $\overline{\partial\Om}=\overline{\Gamma}_1 \bigcup \overline{\Gamma}_2$,
$p(\mathbf{x},t)$ and  $q(\mathbf{x},t)$ are given functions for Dirichlet and Neumann boundary condition respectively.
If the source locations $\mathbf{x}_i^*$ ($i=1,\cdots,s$)
and the corresponding time-dependent intensities $f_i(t)$ ($i=1,\cdots,s$) in (\ref{eq:cfeq}) are all known,
then the distribution of the pollutant concentration $C(\x, t)$ can be obtained by solving the convection-diffusion equation (\ref{eq:cfeq})-(\ref{eq:IC}).
This is usually called a forward or direct problem.
In this paper, we are concerned about the inverse problem, that is,
using the noise-contaminated data $C^{\varepsilon}(\mathbf{x},t)$ ($\varepsilon$ is the noise level) of the concentration $C(\mathbf{x},t)$ in $\Omega$ at terminal time $T$  to recover the source intensity functions $f_i(t)$ ($i=1,\cdots,s$).
In practice, the data $C^{\varepsilon}(\mathbf{x},t)$ is measured by a sensor network placed at some discrete
points inside the domain $\Om$ \cite{L08,NNR98}. A discussion of the sensor network can be found in \cite{L08},
but we shall assume that the measurement data is available here at a set of uniformly
distributed sensors inside $\Om$.

 The Tikhonov optimization algorithm is popular for time-dependent parameter identification problems \cite{KZ98,KZ00}. The main ingredient of the algorithm includes reformulating the reconstruction process as the minimization of  the following functional:
\begin{equation}
J(\mathbf{f})=\frac{1}{2} \int_{\Omega} (C(\mathbf{x}, T)-C^{\varepsilon}(\mathbf{x}, T))^2\,d\x + N_{\beta}(\mathbf{f}),
\label{eq:optim}
\end{equation}
where $\mathbf{f}=(f_1, f_2,\cdots,f_s)^T$, and $N_{\beta}(\mathbf{f})$ denotes some Tikhonov regularization.
Possible choices for the regularizations include $L^2$, $H^1$ and $BV$ regularizations.
Here we consider a combination of  the $L^2$ and  $H^1$ regularizations in the following form
\begin{equation}\label{eq:reguterm}
 N_{\beta}(\mathbf{f}) =\sum_{i=1}^s\frac{\beta_1^i}{2}\int_0^T (f_i(t))^2 dt+ \sum_{i=1}^s\frac{\beta_2^i}{2}\int_0^T |f'_i(t)|^2 dt,
\end{equation}
where $\beta_1^i, \beta_2^i$, $i=1,\cdots,s$, are the $L^2$ or  $H^1$ regularization parameters for the source intensity $f_1(t),\cdots,f_s(t)$ respectively.
The minimization of the functional (\ref{eq:optim}) is subject to the constraints that $C(\x, t)$ satisfies the state equation (\ref{eq:cfeq}) with the boundary conditions (\ref{eq:BCs}) and the initial condition (\ref{eq:IC}).
This has transformed the original inverse source problem into a PDE-constrained optimization problem.
Two kinds of approaches for the optimization problem (\ref{eq:optim}) are available, the discretize-then-optimize approach and the optimize-then-discretize approach.
The solutions from both approaches are credible, although they are not necessarily the same \cite{PBC06}.
We shall use the optimize-then-discretize approach in this work.

Let ${W}^{1,p}(\Om)$ and ${W}^{1,q}(\Om)$ be standard Sobolev spaces with $p, q>0$ such that $1/p+1/q=1$
and $p>2$, $q<2$.
We formally write (\ref{eq:cfeq}) as an operator equation $L(C, \mathbf{f})=0$, then introduce a corresponding Lagrange multiplier $G\in W^{1,p}(\Om)$ and the following Lagrange functional \cite{AB03,KZ98,KZ00}:
\begin{equation}
\mathcal{J}(C,\mathbf{f}, G)=\frac{1}{2}  \int_{\Omega}(C(\x,T)-C^{\varepsilon}(\x, T))^2d\x +N_{\beta}(\mathbf{f})+ (G,L(C,\mathbf{f})),
\label{eq:lagrange1}
\end{equation}
where $G$ is the Lagrange multiplier or the adjoint variable, and $(G,L(C,\mathbf{f}))$ stands for the dual product.

Taking the variations of (\ref{eq:lagrange1}) with respect to $G$, $C$ and $f_i$, $i=1,\cdots,s$,
a system of partial differential equations is derived to characterize the first-order optimality conditions for  this optimization problem (\ref{eq:lagrange1}).
They are the so-called Karush-Kuhn-Tucker (KKT) optimality conditions \cite{KT51}.
It has been verified that the minimization problem (\ref{eq:optim}) is equivalent to solving the KKT system \cite{KT51} of the Lagrangian functional $\mathcal{J}(C,\mathbf{f}, G)$ in \cite{CZ99}.
The three sets of equations in the KKT system are obtained as follows:
\begin{enumerate}
\item[(a)]
The G$\hat{a}$teaux derivative of $\mathcal{J}$ with respect to $G$ at direction $v$ is given by
\begin{align*}
\mathcal{J}_G(C,\mathbf{f}, G)v&=(v,L(C,\mathbf{f}))\\
&=\left(\frac{\partial C}{\partial t},v\right)+(a\nabla C,\nabla v)+(\nabla\cdot (\mathbf{v}C),v)\\
&-\sum_{i=1}^s v(\x_i^*,t)f_i(t)-\langle q,v\rangle_{\Gamma_2}
\end{align*}
for all $v\in W^{1,p}(\Om)$.
\item[(b)]The G$\hat{a}$teaux derivative of $\mathcal{J}$ in (\ref{eq:lagrange1})
with respect to $C$ at direction $w\in W^{1,q}(\Om)$ is given by
\begin{equation}
\begin{split}
\mathcal{J}_{C}(C,\mathbf f, G)w&=\int_{\Omega}(C(\mathbf{x},T)-C^{\varepsilon}(\mathbf{x},T)) wd\x\\
&+ \int_0^T\int_{\Omega} G\left(\frac{\partial w}{\partial t}-\nabla\cdot(a(\mathbf{x})\nabla w)
+\nabla\cdot (\mathbf{v}(\mathbf{x})w)\right) d \x d t.
\end{split}
\label{eq:JwithC}
\end{equation}
For convenience, we write
\[\tilde{L} w:= \frac{\partial w}{\partial t}-\nabla\cdot(a(\mathbf{x})\nabla w)
+\nabla\cdot (\mathbf{v}(\mathbf{x})w,\]
and obtain by integrating by part for the second term of (\ref{eq:JwithC}) that
\begin{align*}
(G,\tilde{L} w)&=\int_0^T\int_{\Omega} G\left(\frac{\partial w}{\partial t}-\nabla\cdot(a(\mathbf{x})\nabla w)+\nabla\cdot (\mathbf{v}(\mathbf{x})w)\right) d \x d t\\
&=\int_{\Omega} G w|_0^T d \x - \int_0^T\int_{\Omega} w\frac{\partial G}{\partial t}d \x d t-\int_0^T\int_{\partial\Omega} a(\mathbf{x})G\frac{\partial w}{\partial \mathbf{n}}d \Gamma d t\\
& +\int_0^T\int_{\Omega} a(\mathbf{x})\nabla w\cdot \nabla G d \x d t+ \int_0^T\int_{\Omega} \nabla\cdot (\mathbf{v}(\mathbf{x})w)G d \x d t\\
&=\int_{\Omega} G w|_0^T d \x+\int_0^T\int_{\partial\Omega} \left(-a(\mathbf{x})\frac{\partial w}{\partial \mathbf{n}}G\right) d \Gamma d t\\
&+\int_0^T\int_{\Omega} \left(-w\frac{\partial G}{\partial t}+a(\mathbf{x})\nabla w\cdot \nabla G + \nabla\cdot (\mathbf{v}(\mathbf{x})w)G\right) d \x d t\,.
\end{align*}
Then applying the boundary and initial conditions of $w$, i.e. $w=0$ on $\Gamma_1$ and $a(\x)\frac{\partial w}{\partial \mathbf{n}}=0$ on $\Gamma_2$, $w(\mathbf{x},0)=0$, we derive
\begin{align*}
(G,\tilde{L}w)&=\int_{\Omega} G(\x,T)w(\x,T) d \x +\int_0^T\int_{\Gamma_1} \left(-a(\mathbf{x})G\frac{\partial w}{\partial \mathbf{n}}\right) d \Gamma d t \\
&+\int_0^T\int_{\Omega} \left(-w\displaystyle\frac{\partial G}{\partial t}+a(\mathbf{x})\nabla w\cdot \nabla G+ \nabla\cdot (\mathbf{v}(\mathbf{x})w)G\right) d \x d t.
\end{align*}
Now noting the arbitrariness of $w$, we can deduce
the adjoint system for the Lagrange multiplier $G$, namely
$G(\x, T)=0$ for $\x\in \Om$, $G(\x, t)=0$ on $\Gamma_1$ and
$G(\x, t)$  satisfies
\begin{equation}
-(G_t, w) + (a\nabla  {G}, \nabla w) + (\nabla\cdot (\mathbf{v}w),G)
=-(\delta(t-T) (C(\cdot,t)-C^{\varepsilon}(\cdot,t)), w)
\label{eq:adjoint1}
\end{equation}
for all  $w\in W^{1,q}(\Om)$ such that $w=0$ on $\Gamma_1$, where
$\delta(t-T)$ is the Dirac delta distribution at $t = T$.
\item[(c)]
The G$\hat{a}$teaux derivative of $\mathcal{J}$ in (\ref{eq:lagrange1}) with respect to $f_i$ at direction $g\in H^1(0,T)$
is given by
\begin{equation}
\begin{split}
{\,\ } \hskip0truecm\label{eq:Jwithf}
\mathcal{J}_{f_i}(C,\mathbf{f}, G)g&=\beta_1^i \int_0^T  f_i(t) g(t) d t + \beta_2^i \int_0^T   f'_i(t) g'(t) d t\\
&-\int_0^T  (G(\x,t), \delta(\x-\x_i^*)g(t)) d t \\
&= \int_0^T (\beta_1^i f_i(t) -  G(\x_i^*,t))g(t) d t +  \beta_2^i \int_0^T  f'_i(t)g'(t) d t.
\end{split}
\end{equation}
\end{enumerate}
Putting (a)-(c) together, the KKT system is formulated as follows:
\begin{align}
\begin{cases}
\mathcal{J}_G(C,\mathbf{f},G)v = 0\\
\mathcal{J}_C(C,\mathbf{f}, G)w = 0\\
\mathcal{J}_{f_i}(C,\mathbf{f}, G)g = 0,\quad i=1,\cdots,s,
\end{cases}
\end{align}
that is, for any $v\in W^{1,p}(\Om)$ and $w\in W^{1,q}(\Om)$, we have the following coupled system:
\begin{align}
\begin{cases}
\left(\displaystyle\frac{\partial C}{\partial t},v\right)+(a\nabla C,\nabla v)+(\nabla\cdot (\mathbf{v}C),v)-\displaystyle\sum_{i=1}^s v(\x_i^*)f_i(t)-\langle q,v\rangle_{\Gamma_2} = 0\\
-\left(\displaystyle\frac{\partial G}{\partial t},w\right)+(a(\mathbf{x})\nabla G, \nabla w)+(\nabla\cdot(\mathbf{v}(\x)w), G)\\
+(\delta(t-T) (C(\cdot,t)-C^{\varepsilon}(\cdot,t), w)) = 0\\
- ( G(\x_i^*,\cdot), g) +  \beta_1^i (f_i, g) +  \beta_2^i ( f'_i, g')=0, \quad i=1,\cdots,s
\end{cases}
\label{eq:kkt}
\end{align}
with $C(\x, 0)=C_0(\x)$, $G(\x, T)=0$.
The rest of the paper is devoted to solving (\ref{eq:kkt}) as a coupled space-time system.
It is noted that the first equation in (\ref{eq:kkt}) is the state equation,
and the second equation is the adjoint equation,
and the last set of equations are elliptic equations with respect to each unknown source intensity.

\section{Finite element discretization}\label{sec:finite element}
Let $\mathcal{T}^h$ be a triangulation of $\Omega$ with triangular elements,
then we define $V^h$ as the finite element space \cite{C78} consisting of continuous piecewise linear functions on $\mathcal{T}^h$, and $\mathring{V}^h$ the subspace of $V^h$ with functions vanishing on the Dirichlet boundary $\Gamma_1$.
To fully discretize the system (\ref{eq:kkt}), we partition the time interval $[0,T]$ as $0=t^0<t^1<\cdots< t^M=T,$ with $t^n=n \tau, \tau=T/M$.
Define $U^{\tau}$ as a piecewise linear continuous finite element space in time.
For a given sequence $\{H^n(\x) = H(\x, t^n)\}$, we define the difference quotient and  the averaging function respectively by
\begin{equation}\label{eq:timepart}
\partial_{\tau} H^n (\x)= \frac{H^n(\x)-H^{n-1}(\x)}{\tau},\qquad \bar{H}^n=\frac{H(\x, t^{n-1})+H(\x, t^n)}{2}.
\end{equation}
Let $\pi_h$ be the finite element interpolation associated with the space $V^h$, and $C_h^n(\x)$ be the finite element approximation of $C(\x, t^n)$,
then we discretize the state and adjoint equations of the system (\ref{eq:kkt}) by the Crank-Nicolson scheme in time and piecewise linear finite elements  in space,
and lastly we use piecewise linear finite element in time to discretize the equations with respect to $\mathbf{f}$.
The finite element approximation of the KKT system (\ref{eq:kkt}) can be formulated as follows:

Find a sequence of approximations $C_h^n$, $G_h^n\in V^h$ ($0\le n\le M$), $f_i^{\tau}\in U^{\tau}, i=1,\cdots,s$, such that $C_h^0=\pi_hC_0$, $G_h^M=\mathbf{0}$ and $C_h^n(\x)=\pi_h p(\x, t^n), G_h^n(\x)= \mathbf{0}$ for $\x\in \Gamma_1$ satisfying
\begin{align}
\begin{cases}
(\partial_{\tau} C_h^n, v_h) + (a\nabla  \bar{C}_h^n, \nabla v_h) + (\nabla\cdot (\mathbf{v}\bar{C}_h^n),v_h) \\
= \sum_{i=1}^s v_h(\x_i^*) \bar{f}_i^n+\langle\bar{q}^n,v_h\rangle_{\Gamma_2}, ~~\forall\,v_h\in \mathring{V}^h,~ n=1,\cdots,M\\
-(\partial_{\tau} G_h^n, w_h) + (a\nabla  \bar{G}_h^n, \nabla w_h) + (\nabla\cdot (\mathbf{v}w_h),\bar{G}_h^n)\\
=-\chi_n ((C_h^n-C^{\varepsilon}), w_h),~~\forall\,w_h\in \mathring{V}^h,~ n=M,\cdots,1\\
-  (G^{\tau}(\x_i^*, \cdot), g^{n}) +\beta_1^i( f^{\tau}_i,g^{n}) +  \beta_2^i ((f^{\tau}_i)', (g^{n})')=0,~~ i=1,\cdots,s, ~ n=0,\cdots,M,
\end{cases}
\label{eq:kkt2}
\end{align}
where $\chi_n = \frac{1}{2}$ when $n=M$ and 0 when $1 \leq n< M$.
We denote the basis functions of finite element spaces $V^h$ and $U^{\tau}$ by $\phi_i$, $i=1,\cdots,N$ and $g^n$, $n=0,\cdots,M$, respectively,
and introduce the following matrices:
\begin{align*}
A&=(a_{ij})_{i,j=1,\cdots,N},\qquad\quad a_{ij}=(a \nabla \phi_i, \nabla \phi_j)\\
B&=(b_{ij})_{i,j=1,\cdots,N},\qquad\quad b_{ij}=( \phi_i,  \phi_j)\\
E&=(e_{ij})_{i,j=1,\cdots,N},\qquad\quad e_{ij}=(\nabla\cdot(\mathbf{v} \phi_i), \phi_j)\\
K&=(k_{nm})_{n,m=0,\cdots,M},\qquad k_{nm}=( (g^n)',  (g^{m})')\\
D&= (d_{nm})_{n,m=0,\cdots,M},\qquad d_{nm}=( g^n, g^{m})\\
\end{align*}
and the vectors
\begin{align*}
C^n&=(C_1^n,C_2^n,\cdots,C_N^n)^T,\qquad \mbox{for}~~n=0,\cdots,M\\
G^n&=(G_1^n,G_2^n,\cdots,G_N^n)^T,\qquad \mbox{for}~~n=0,\cdots,M\\
f_k&=(f_k^0, f_k^1,\cdots,f_k^M)^T, \qquad \mbox{for}~~ k=1,\cdots,s\\
r^n&=(r_1^n, r_2^n,\cdots,r_N^n)^T, \quad \mbox{for}~~ j=1,\cdots,N,~~ n=1,\cdots,M~~ \mbox{with}\\
r_j^n&=-\tau\left(\sum_{k=1}^s \phi_j(\x_k^*) \frac{(f_k^{n}+f_k^{n-1})}{2}+\left\langle \frac{(q^{n}+q^{n-1})}{2},\phi_j\right\rangle_{\Gamma_2}\right)\\
g_k^*&=(g_k^0, g_k^1,\cdots,g_k^M)^T,~~ \mbox{with}~ g_k^n=G(\x_k^*, t^n), \qquad \mbox{for}~~ k=1,\cdots, s,~~n=0,\cdots,M \\
d&=(d_1, d_2,\cdots,d_N),~~ \mbox{with}~d_j= (C^{\varepsilon}, \phi_j) \qquad \mbox{for}~~ j=1,\cdots,N.
\end{align*}
The matrix form of the KKT system is then reformulated by using the above notations as the following:
\begin{align}
\begin{cases}
\left(B+\displaystyle\frac{\tau}{2}(A+E)\right)C^{n}+\left(-B+\displaystyle\frac{\tau}{2}(A+E)\right)C^{n-1}+r^{n}=0,  ~~ n=1,\cdots,M \\
\left(-B+\displaystyle\frac{\tau}{2}(A+E^T)\right)G^{n}+\left(B+\displaystyle\frac{\tau}{2}(A+E^T)\right)G^{n-1}\\
+ \tau\chi_n (B C^{n}-d)=0,  ~~ n=M,\cdots,1\\
-D g_k^*+(\beta_1^k D+ \beta_2^k K) f_k =0,  ~~ k=1,\cdots,s.
\end{cases}
\label{eq:kkt4}
\end{align}
We can follow the approaches in \cite{JZ10,KZ98,KZ00,XZ05} to obtain the convergence of the discretized problem (\ref{eq:kkt2}) to the continuous optimization problem (\ref{eq:lagrange1}).
\section{A space-time domain decomposition method for the KKT system}\label{sec:scalable solvers}
\subsection{Fully coupled KKT system with special ordering of unknowns}
The ordering of the unknowns for the discretized KKT system (\ref{eq:kkt2}) has significant influence in the convergence and computing efficiency of the iterative solver.
Traditional reduced space SQP methods split the system into three subsystems and solve each subsystem for $C$, $G$, and $\mathbf{f}$ one by one in sequential steps \cite{DZZ12},
in this case the unknowns are ordered physical variable by physical variable.
To develop a scalable and fully coupled method for solving the KKT system,
we use the so-called fully coupled ordering, the unknowns $C$ and $G$ are ordered mesh point by
mesh point and time step by time step. At each mesh point $\x_{j}, j=1,\cdots,N$, and time step $t^n, n=0,\cdots,M$, the unknowns are ordered in the order of $C_{j}^n, G_{j}^n, j=1,\cdots,N, n=0,\cdots,M$.
Such ordering contains unknowns of the same space-time subdomain in a subblock, preconditioners such as additive Schwarz can be applied naturally to each subblock of the fully coupled KKT system and the ordering also improves the cache performance of the LU factorization based solvers.
Since $\mathbf{f}$ is defined only in the time dimension, we put all the unknowns of $\mathbf{f}$ at the end after $C$ and $G$. More precisely, we define the solution vector $U$ by
\begin{align*}
U&=(C_{1}^0,G_{1}^0,\cdots,C_{N}^0, G_{N}^0, C_{1}^1, G_{1}^1,\cdots,C_{N}^1, G_{N}^1,\cdots,C_{1}^M,G_{1}^M, \\
 & \cdots,C_{N}^M, G_{N}^M, f_1^0,\cdots,f_s^0,\cdots,f_1^M,\cdots,f_s^M)^T
\end{align*}
then the linear system (\ref{eq:kkt2}) with unknowns $C_j^n$ and $G_j^n$, $j=1,\cdots,N$, $n=0,\cdots,M$, and $f_k^n$, $k=1,\cdots,s$, $n=0,\cdots,M$, is reformulated into the following linear system:
\begin{equation}\label{eq:kktsys}
\mathbf{F}U=b,
\end{equation}
where $\mathbf{F}$ is a sparse matrix of size $(M+1)(2N+s) \times (M+1)(2N+s)$ derived from the finite element discretization for KKT system (\ref{eq:kkt4}) with the following block structure:
\[\mathbf{F}=\left(\begin{array}{cccccc}S_{00}&S_{01}&\mathbf{0}&\cdots&\mathbf{0}&S_{0,M+1}\\ S_{10}&S_{11}&S_{12}&\cdots&\mathbf{0}&S_{1,M+1}\\ \vdots& \ddots & \ddots&\ddots &\vdots&\vdots \\
\mathbf{0}&\cdots&S_{M-1,M-2}&S_{M-1,M-1}&S_{M-1,M}&S_{M-1,M+1}\\ \mathbf{0}&\cdots&\mathbf{0}&S_{M,M-1}&S_{M,M}&S_{M,M+1}\\ S_{M+1,0} & S_{M+1,1}&S_{M+1,2}&\cdots&S_{M+1,M}& S_{M+1,M+1}\end{array}\right),\]
and $b$ has the following form correspondingly:
\[ b=(b_{0}, b_{1},\cdots,b_{M+1})^T.\]
In the matrix $\mathbf{F}$, the block matrices $S_{ij}$, with $0\leq i,j \leq M$ are of size $2N\times2N$ and are zero matrices except the ones in tridiagonal stripes $\{S_{i,i-1}\}, \{S_{i,i}\}, \{S_{i,i+1}\}$.
The stripe $\{S_{i, M+1}\}$, $0\leq i\leq M$ are nonzero sparse blocks of size $2N\times s(M+1)$; furthermore $\{S_{M+1,i}\}$, $0\leq i\leq M$ are nonzero sparse blocks of size $s(M+1)\times 2N$ and $S_{M+1,M+1}$ is a nonzero
tridiagonal matrix of size $s(M+1)\times s(M+1)$.

\subsection{Space-time Schwarz preconditioners}\label{sec:space-time schwarz}
The KKT system (\ref{eq:kktsys}) is usually large in size and severely ill-conditioned.
In traditional reduced space SQP methods, the subsystems corresponding to unknowns $C$ and $G$ are time-dependent,
time-marching algorithms starting from the initial or terminal moment are applied.
But we notice that
the adjoint equation needs the concentration distribution of $C$ at the terminal time  $t=T$, which means that the state equation and the adjoint equation should be solved in a sequential order.
In addition, sequential steps within reduced space SQP methods exist between both the KKT subsystems and the time marching for time-dependent inverse problems, thus are quite challenging for efficient parallelization.
To overcome the lack of parallelism in SQP methods, we shall propose to solve the fully coupled system (\ref{eq:kktsys}) all at once.
This is a very large system, the all-at-once method is traditionally regarded as a very expensive approach and not suitable for small computers.
But on high-performance computers, especially on the upcoming exascale computers,
we believe this approach is more attractive than the reduced space methods.
It is well known that a direct solver such as Gaussian elimination or LU factorization is not suitable for very large problems
due to the lack of parallel scalability.
We shall use a preconditioned Krylov subspace method, where some preconditioning technique
will be introduced for reducing the condition number of the KKT system and accelerating the convergence rate of the Krylov
subspace method.
Various preconditioners have been developed and applied for various elliptic and parabolic systems,
such as the (block) Jacobi method, (incomplete) LU factorization, (multiplicative) additive Schwarz method, multigrid method, multilevel method, etc. \cite{B96,BG99,CLZ09}.
Among these preconditioners  the Schwarz type domain decomposition method is shown to have excellent preconditioning effect and parallel scalability \cite{CLZ09,PBC06}.

We shall propose a ``space-time'' Schwarz type preconditioner for the unsteady inverse problems.
Different from the classical Schwarz type preconditioning technique which only decomposes the space domain,
we want full parallelization in both space and time.
The idea of space-time parallel algorithm comes from the parareal algorithm, proposed by Lions et al. in \cite{LMT01}.
The parareal algorithm is an iterative method which involves a coarse (coarse grid in the time dimension) solver for prediction and a fine (fine grid in the time dimension) solver for correction.
An insight on the stability and convergence of the parareal algorithm was given in \cite{GH08,SR05}.
Parareal algorithm has been applied to solve problems in  molecular dynamics \cite{BBMTZ02}, fluid and structural mechanics \cite{FC03}, quantum control \cite{MT03} etc.
However in the implementation of the parareal algorithm, the scalability is determined largely by the coarse time step and the space discretization scheme.
In \cite{MT05}, the parareal algorithm was combined with domain decomposition in space to achieve higher degree of parallelization.
Different from the parareal algorithm, the new ``space-time'' Schwarz type preconditioner treats the time variable and the space variables equally, so the physical domain is a ``space-time'' domain, instead of the conventional space domain. We apply a domain decomposition technique to the coupled ``space-time'' domain.

In each ``space-time'' subdomain, a time-dependent subproblem with vanishing space boundary conditions
and vanishing data at ``artificial'' initial and terminal time is solved.
The same as the global problem, no time-marching is performed in each subproblem, all unknowns associated to the same
space-time subdomain are solved simultaneously.
The proposed ``space-time'' Schwarz preconditioner eliminates all sequential steps and all unknowns are treated at the same level of priority.
We use a right-preconditioned restarted GMRES to solve the system (\ref{eq:kktsys}):
 \[\mathbf{F}M^{-1}U'=b,\]
where $M^{-1}$ is a ``space-time'' additive Schwarz preconditioner and $U=M^{-1}U'$.

\begin{figure}
\centering
\includegraphics[width=8cm]{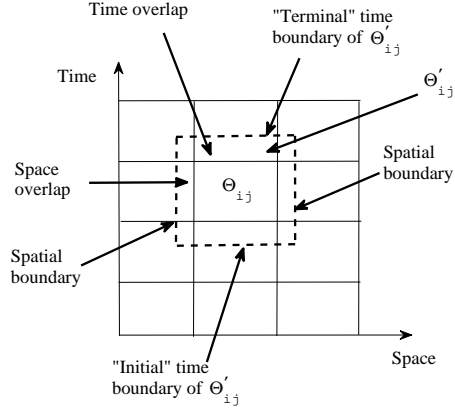}
\vskip-0.6truecm
\caption{``Space-time'' domain decomposition - an overlapping subdomain with boundary conditions.}\label{fig:ddm}
\end{figure}

To formally define the preconditioner $M^{-1}$ we need to introduce a partition of the space-time domain
$\Om\times[0,T]$, denoted by $\Theta$.
Firstly we decompose the domain $\Om$ into nonoverlapping subdomains $\Om_i$, $i=1,\cdots,N_1$, and then divide the time interval $[0,T]$ into subintervals $T_j=[t_{j-1}, t_j]$, $j=1,2,\cdots,N_2$, and $0=t_0<t_1<\cdots<t_{N_2}=T$.
We remark that the time partition here is coarser than that used in the full finite element discretization
described in Section~\ref{sec:finite element},
and each interval $T_j$ contains a few consequent time intervals $[t^k,t^{k+1}]$.
$\Theta$ consists of $\Theta_{ij}=\Om_i\times T_j$, $i=1,\cdots,N_1, j=1,\cdots,N_2$.
In order to obtain an overlapping decomposition of $\Theta$, we extend each subdomain $\Om_i$ to a larger region $\Om_i'$ and each subinterval $T_j$ to a longer interval $T_j'$, satisfying $\Om_i \subset \Om_i'$, $T_j \subset T_j'$. Now each $\Theta_{ij}$ can be straightforwardly extended to $\Theta_{ij}' = \Om_i'\times T_j'$ with $\Theta_{ij} \subset \Theta_{ij}'$.
The sizes of $\Theta_{ij}'$  are chosen so that the overlap is as uniform as possible around the perimeter of interior domains $\Theta_{ij}' \subset \Theta$.
For boundary subdomains we neglect the part outside of $\Theta$.
See Figure \ref{fig:ddm} for an illustration of the space-time domain decomposition.

We denote the size of the KKT matrix $\mathbf{F}$ by $\tilde{N}\times \tilde{N}$, clearly there are two degrees of freedom at each mesh point corresponding to the state variable $C$ and the adjoint variable $G$.
The unknown time-dependent source intensity variables are allocated on
the same processor as  the last space-time subdomain $\Theta'_{N_1,N_2}$.

 On each extended subdomain $\Theta_{ij}'$, we define the $\tilde{N}_{ij}\times \tilde{N}$ matrix $R^{\delta}_{ij}$, its $2\times 2$ block element $(R^{\delta}_{ij})_{l_1, l_2}$ is either an identity block if the integer indices $l_1$ and $l_2$ are related to the same mesh point and time step and they belong to $\Theta_{ij}'$ or a zero block otherwise.
 The multiplication of $R^{\delta}_{ij}$ with an $\tilde{N}\times 1$ vector generates a shorter vector by keeping all components corresponding to the subdomain $\Theta_{ij}'$.
 $R^{0}_{ij}$ is defined similarly as $R^{\delta}_{ij}$, with the difference that its application to a $\tilde{N}\times 1$ vector excludes the mesh points in $\Theta_{ij}' \backslash \Theta_{ij}$.
Now for each space-time subdomain we have defined the following local problem:
\begin{equation}
\begin{cases}
\displaystyle\frac{\partial G}{\partial t}=-\nabla \cdot (a(\mathbf{x})\nabla G)-\mathbf{v}(\mathbf{x})\cdot \nabla G \\
{\,\ } \hskip0truecm +\delta(t-T)(C(\x,t)-C^{\varepsilon}(\x,t)), ~ ~ (\x,t)\in \Theta_{ij}' \\
\displaystyle\frac{\partial C}{\partial t}=\nabla\cdot(a(\mathbf{x})\nabla C)-\nabla\cdot (\mathbf{v}(\mathbf{x})C)+\sum_{i=1}^s \delta(\mathbf{x}-\mathbf{x}_i^*)f_i(t), ~ (\x,t)\in \Theta_{ij}'.
\end{cases}
\label{eq:subprob}
\end{equation}
It is complemented by the following boundary conditions
\begin{align}\label{eq:subprobBCs}
C(\mathbf{x},t)=0; \q G(\mathbf{x},t)=0,~ \x\in \partial\Om_i'
\end{align}
along with the ``initial'' and ``terminal'' time boundary conditions
\begin{equation}\label{eq:subprobIC}
C(\mathbf{x},t_{j-1})=0; \q G(\x, t_{j-1})=0
\end{equation}
\begin{equation}\label{eq:subprobTC}
C(\mathbf{x},t_{j})=0; \q G(\x, t_{j})=0.
\end{equation}
For the last subdomain, we include the additional variables corresponding to the source intensities $f_i$, $i=1,\cdots,s$ satisfying
\begin{equation}\label{eq:subprobf}
\beta_2^i  f''_i +\beta_1^i f_i + G(\x^*, \cdot) = 0,
\end{equation}
with the Neumann condition
\begin{equation}\label{eq:subprobfIC}
f'_i(t)=0, \qquad t=0, T.
\end{equation}
We remark that (\ref{eq:subprob}) is a parabolic system and it is usually ``illegal'' to impose both initial and terminal conditions (\ref{eq:subprobIC})-(\ref{eq:subprobTC}). However, as inexact local solvers on space-time subdomains that form
the global preconditioner, such local boundary conditions work well as we shall see from
our numerical experiments. Similar boundary conditions are used in the context of hyperbolic subdomain problems
\cite{YCC10}.

 Let $M_{ij}$ be a discretization of (\ref{eq:subprob})-(\ref{eq:subprobTC}) and $M_{ij}^{-1}$ be an exact or approximate inverse of $M_{ij}$.
 The space-time additive Schwarz preconditioner is defined as
 \[{\,\ } \hskip0truecm M_{asm}^{-1} = \sum_{j=1}^{N_2}\sum_{i=1}^{N_1} (R^{\delta}_{ij})^T M_{ij}^{-1}R^{\delta}_{ij}. \]
 It is noted that the last space-time subdomain solver $M_{N_1,N_2}^{-1}$ is an inverse or an approximate inverse of the matrix arising from the discretization of  the subproblem (\ref{eq:subprob})-(\ref{eq:subprobTC}) of $\Theta'_{N_1, N_2}$  and (\ref{eq:subprobf})-(\ref{eq:subprobfIC}).
 Although its construction is slightly different from that of the other subdomain inverse matrices, we still use the same notation.

In addition to the standard additive Schwarz method (ASM) described above,
the restricted version (RAS) of the method developed in \cite{CS99} for standard space domain decompositions is also widely used.
So we extend it to our current space-time domain decomposition,
then the space-time RAS preconditioner is defined as
\[{\,\ } \hskip0truecm M_{ras}^{-1} = \sum_{j=1}^{N_2}\sum_{i=1}^{N_1} (R^{\delta}_{ij})^TM_{ij}^{-1}R^{0}_{ij}.\]
For some applications, RAS achieves better preconditioning effect with less communication time since one of the restriction or extension operations does not involve any overlap. We use the restricted version in our experiments to be presented in the next section.

We remark that, computationally, the matrix $M_{ij}$ can be obtained as $R^{\delta}_{ij}\mathbf{F}(R^{\delta}_{ij})^T$.
Moreover, if $N_2=1$, then no time partition is performed in the time dimension, $M_{ras}^{-1}$ is a ``space-only'' domain decomposition preconditioner for the fully coupled KKT system.

\section{Numerical examples}\label{sec:numerical examples}
We present in this section some numerical examples of recovering the intensity functions $f_i(t)$ ($i=1,\cdots,s$) at given source locations $\x_1^*,\cdots,\x_s^*$.
We set the test domain to be $\Om=(-2, 2)\times (-2, 2)$, and the terminal time at $T=1$.
We denote the time step by $n_t$ and the number of mesh points in $x$ and $y$ directions
by $n_x$ and $n_y$, respectively.
Homogeneous Dirichlet and Neumann boundary conditions are imposed
on $\Gamma_1=\{\mathbf{x}=(x_1,x_2); ~|x_1|=2\}$  and  $\Gamma_2=\{\mathbf{x}=(x_1,x_2);  ~|x_2|=2\}$, respectively.
The diffusive coefficient $a(\mathbf{x})$ and the convective coefficient $\mathbf{v}(\x)$ are chosen to be $1.0$ and $(1.0,1.0)^T$, respectively.

The preconditioned KKT coupled system will be solved by the restarted GMRES method with a restart number 50 \cite{S03}.
For clarity, we provide the restarted GMRES algorithm below for a general linear system $Ax=b$:

\begin{algorithm}[H]
\caption{ Restarted GMRES method with a restart number $m$}
\label{alg:rgmres}
\begin{algorithmic}
      \STATE      1.~Compute $r_0=b-Ax_0$, $\beta=\|r_0\|_2$, and $v_1=r_0/\beta$;
      \STATE      2.~Generate the Arnoldi basis and the matrix $\bar{H}_m$ using the Arnoldi algorithm
      \STATE  \quad starting with $v_1$;
      \STATE      3.~Compute $y_m$ which minimizes $\|\beta e_1-\bar{H}_m y\|_2$ and $x_m=x_0+V_m y_m$;
       \STATE     4.~Stop the iteration if the stopping criteria are satisfied; otherwise set $x_0:=x_m$
       \STATE \quad and go to 1.
\end{algorithmic}
\end{algorithm}
\noindent
In the algorithm above, $e_1$ is the first column of the identity matrix, $V_m$ stands for
the $n\times m$ matrix with columns $v_1, \cdots, v_m$, and $\bar{H}_m$ denotes the $(m+1)\times m$ Hessenberg matrix
\cite{S03}.
The computational cost is overwhelming and becomes more unstable numerically
when $m$ is large. So we should set $m$ to a reasonable number and
get restarted when the stopping criteria are not satisfied.

For definiteness, we shall denote as {\it a cell}
the smallest space-time element after the space triangulation of $\Om$ and time partition of the interval $[0, T]$.
The size of overlap, that is the number of overlapping cells, is denoted by $iovlp$ and set to $4$ unless otherwise specified.
The subsystem is solved with a sparse LU factorization or an incomplete LU factorization (ILU) with the fill-in level denoted by $ilulevel$.
We still take the linear system $Ax=b$ for example to explain the definition of the fill-in level of an ILU factorization.
Based on the Gaussian elimination, each location $(i,j)$ of the sparse matrix $A$ has a level of fill,
denoted by $\mbox{lev}_{ij}$, which should indicate that the higher the level is, the smaller the element is
\cite{S03}:
\begin{enumerate}
\item The initial value of fill of an element $a_{ij}$ of $A$ is defined by
\[\mbox{lev}_{ij}=\left\{\begin{array}{cc}0 & \mbox{if}~ a_{ij}\neq 0~\mbox{or}~i=j\\ \infty &\mbox{otherwise}\end{array}\right.\]
\item Each time in the process of Gaussian elimination, the element $a_{ij}$ is updated in the loop of $k$ by
$a_{ij}=a_{ij}-a_{ik}a_{kj}$,
the corresponding level of fill is updated by
$\mbox{lev}_{ij}=\mbox{min}\{\mbox{lev}_{ij}, \mbox{lev}_{ik}+\mbox{lev}_{kj}+1\}.$
\end{enumerate}
In ILU with fill-in level $p$, an element whose level of fill $\mbox{lev}_{ij}$ does not exceed $p$ will be kept.
So the larger the level of fill $p$ is, the more elements are kept in the factorization.

For the scalability test we use LU factorization as the subdomain solver and incomplete LU factorization with $ilulevel=3$ for the other tests if not specified.
The relative residual convergence tolerance of GMRES is set to be $10^{-5}$.
The algorithm is implemented based on the package Portable, Extensible Toolkit for
Scientific computation (PETSc) \cite{BBEG10} and run on a Dawning TC3600 blade server system at the National Supercomputing Center in Shenzhen,
China with a 1.271 PFlops/s peak performance.

We use a high resolution numerical solution of the concentration at the terminal time $t=T$ as the noise-free observation data.
In other words, we first solve the forward convection-diffusion system (\ref{eq:cfeq})-(\ref{eq:IC}) on a very fine mesh,  $640\times 640$,  with a small time stepsize $\tau=1/160$, then add a random noise of the following form to the terminal concentration
\[
{\,\ } \hskip0truecm C^{\varepsilon}(\mathbf{x})=(1+\varepsilon\, r)C(\mathbf{x},T),
\]
where $r$ is a random function with uniform distribution in $[-1,1]$, and $\varepsilon$ is the noise level.
In our numerical experiments, $\varepsilon$ is set to 1\% if not specified.
\subsection{Reconstruction results and parallel efficiency tests}
The tests are designed to investigate the recovery effect of the pollutant source intensity functions and
to understand how the solution of the KKT system behaves when using different mesh sizes, time steps,
regularization parameters and number of processors, which is denoted by $np$.
Supposing the source locations are known, we consider the following four examples.
\begin{enumerate}
\item[(1)]  $f=t^2, \quad\qquad\qquad\qquad\qquad\quad\quad\quad\quad\,\,\, \x_1^*=(1.0, 1.0)^T $.
\item[(2)]  $f=\displaystyle\frac{75}{4}t(1-t)\left(\displaystyle\frac{1}{6}-t\right)^2+1.0, \qquad\quad \x_1^*=(1.0, 1.0)^T $.
\item[(3)]
$f_1=t^2, \quad\qquad\qquad\qquad\qquad\qquad\qquad\, \x_1^*=(1.0, -1.0)^T$\\
$f_2=\displaystyle\frac{75}{4}t(1-t)\left(\displaystyle\frac{1}{6}-t\right)^2+1.0, \quad\quad\,\,\, \x_2^* =(0.0, 0.0)^T$.
\item[(4)]
$f_1=t^2, \quad\qquad\qquad\qquad\qquad\qquad\qquad\, \x_1^*=(34/79, 24/79)^T$\\
$f_2=\displaystyle\frac{75}{4}t(1-t)\left(\displaystyle\frac{1}{6}-t\right)^2+1.0, \quad\quad\,\,\, \x_2^* =(14/79, 14/79)^T$\\
$f_3=3-t, \quad\qquad\qquad\qquad\qquad\quad\quad\,\,\,\,\, \x_3^*=(25/79, 15/79)^T$.\\
\end{enumerate}

\begin{figure}
\centering
\includegraphics[width=8cm]{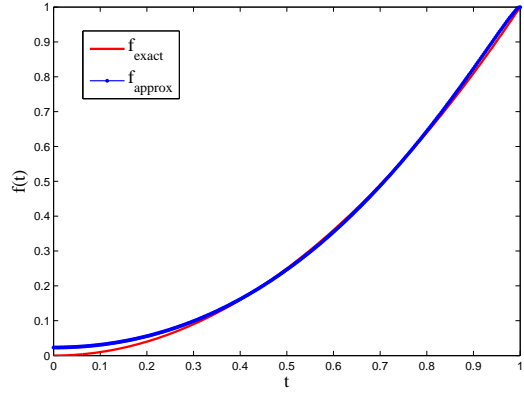}
\vspace{-0.15in}
\caption{Comparison of analytical and computed solution for Example 1.}
\label{fig:ex1}
\end{figure}

\textbf{Example 1}.
This is an example of recovering a quadratic polynomial source intensity function.
An $H^1$ regularization is applied and the parameter is chosen to be $\beta_2=10^{-4}$ ($\beta_1=0$).
Figure \ref{fig:ex1} shows the reconstructed result with mesh $n_x=80, n_y=80$ and time step $n_t=320$, when 64 processors are used.
The blue dotted line represents the reconstructed source intensity which is quite close to the red true shape.
This shows that the time-dependent intensity is successfully recovered by the algorithm.

\begin{figure}
\centering
\includegraphics[width=6cm]{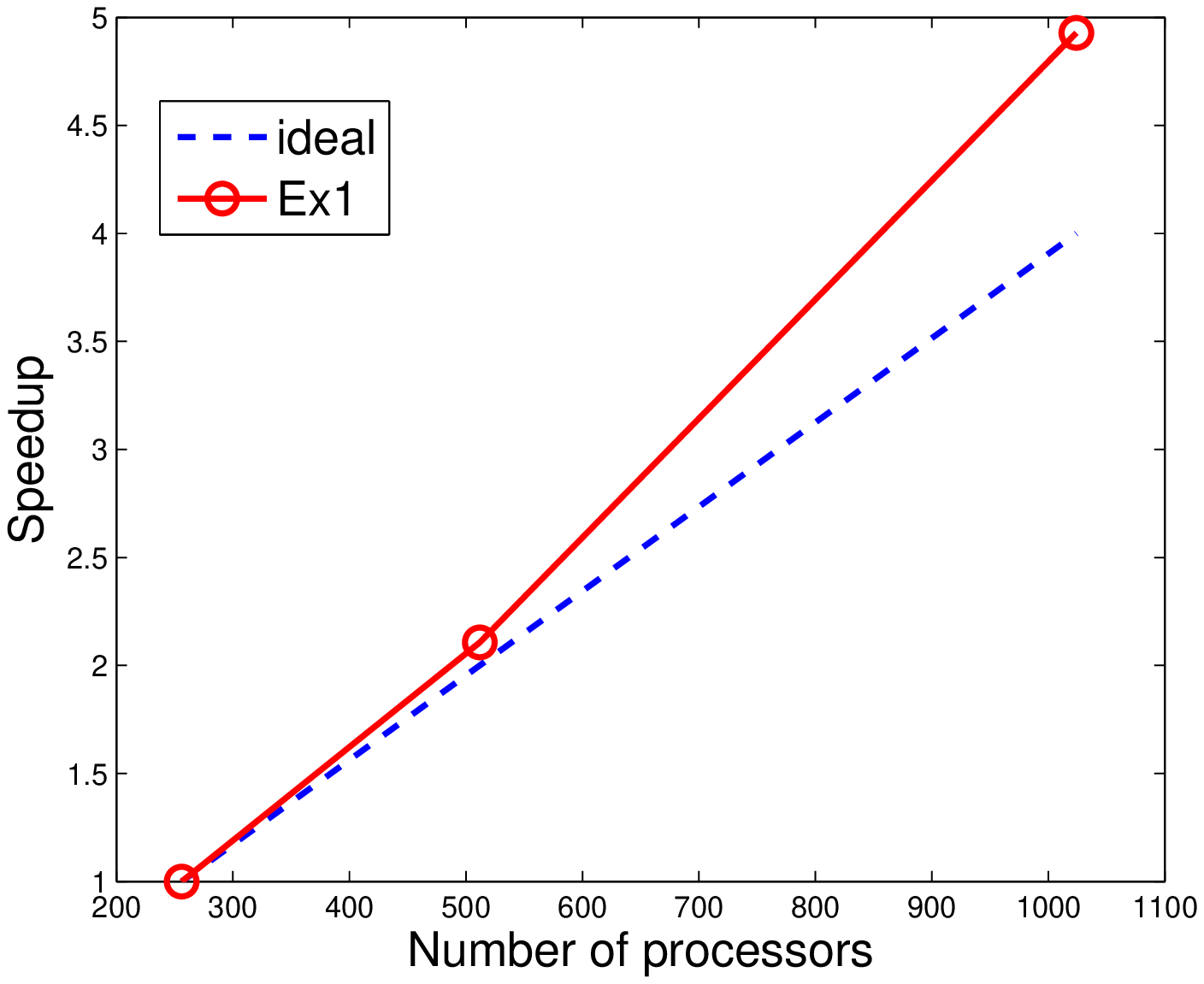}
\includegraphics[width=6cm]{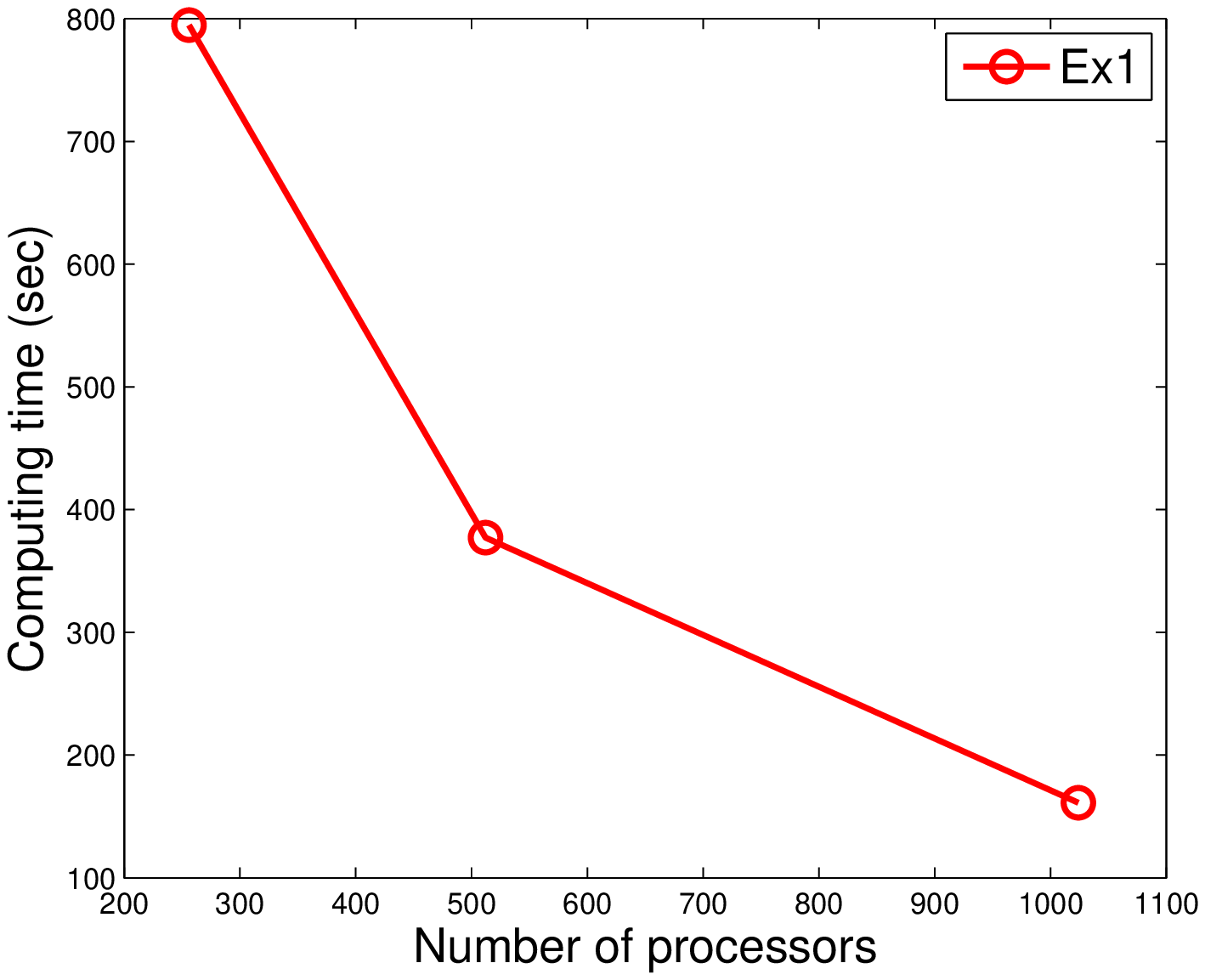}
\vspace{-0.15in}
\caption{The speedup (left) and computing time (right) for Example 1.}
\label{fig:scalex1}
\end{figure}

We present the strong scalability results in Table \ref{tab:mesh1}. Sparse LU factorization is applied as the subdomain solver.
The spatial mesh is $160\times160$ and the number of time steps is 160. The total degrees of freedom is $8,192,160$.
As the number of processors increases, the computing time decreases significantly and superlinear speedup is obtained, for $np\leq 1024$, in Figure \ref{fig:scalex1}.
Since the number of processors is the same as the number of subdomains, more processors lead to an increasing number of iterations.
This suggests that the condition number of the preconditioned KKT matrix depends on the number of subdomains. Similar dependency was proved for elliptic problems \cite{TW04}.

We fix the number of processors to $np=128$, the mesh to $n_x=80$, $n_y=80$ and the time step $n_t=320$, then test several choices of regularization parameters.
ILU factorization is used as the subdomain solver with the fill-in level being $ilulevel=3$.
From the results in Table \ref{tab:ex1meshregu}, as $\beta_2$ becomes smaller, the number of GMRES iterations increases, and no significant change is observed for the total computing time.

\begin{table}
\centering
\begin{tabular}{|ccccc|}
\hline
$np$&Its&Time(sec)&Speedup&Ideal\\
\hline
256&145&794.70&1&1\\
512&181&377.16&2.11&2\\
1024&242&161.23&4.93&4\\
\hline
\end{tabular}
\vspace{0.1in}
\caption{Scalability test for Example 1: $n_t=160,~ n_x=160,~ n_y=160,~ DOF=8,192,160$.}
\label{tab:mesh1}
\end{table}

\begin{figure}
 \centering
\includegraphics[width=8cm]{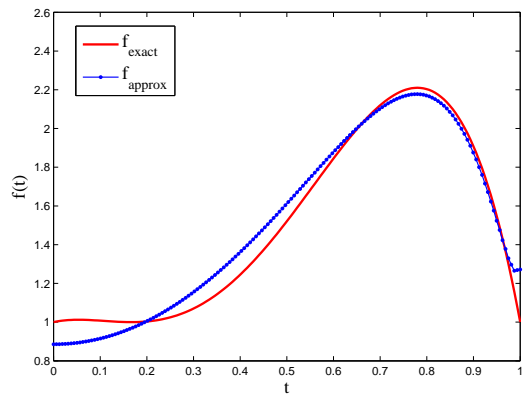}
\vspace{-0.15in}
\caption{Comparison of analytical and computed solution for Example 2.}
\label{fig:ex2}
\end{figure}

\begin{table}
\centering
\begin{tabular}{|ccc|}
\hline
$\beta_2$&Its&Time(sec)\\
\hline
$10^{-4}$&88&54.91\\
$10^{-5}$&93&55.75\\
$10^{-6}$&96&56.68\\
\hline
\end{tabular}
\vspace{0.1in}
\caption{$H^1$ regularization parameter test for Example 1: $n_t=320, n_x=80, n_y=80, DOF=4,096,320, np=128$.}
\label{tab:ex1meshregu}
\end{table}
\textbf{Example 2}.
This is an example of recovering a polynomial source intensity function of degree 4.
We set $\beta_1=0$ and use an $H^1$ regularization with $\beta_2=10^{-4}$.
Satisfactory result is shown in  Figure \ref{fig:ex2} with mesh $n_x=80, n_y= 80$ and time step $n_t=160$,
when 64 processors are used for the computation.

\begin{figure}
\centering
\includegraphics[width=6cm]{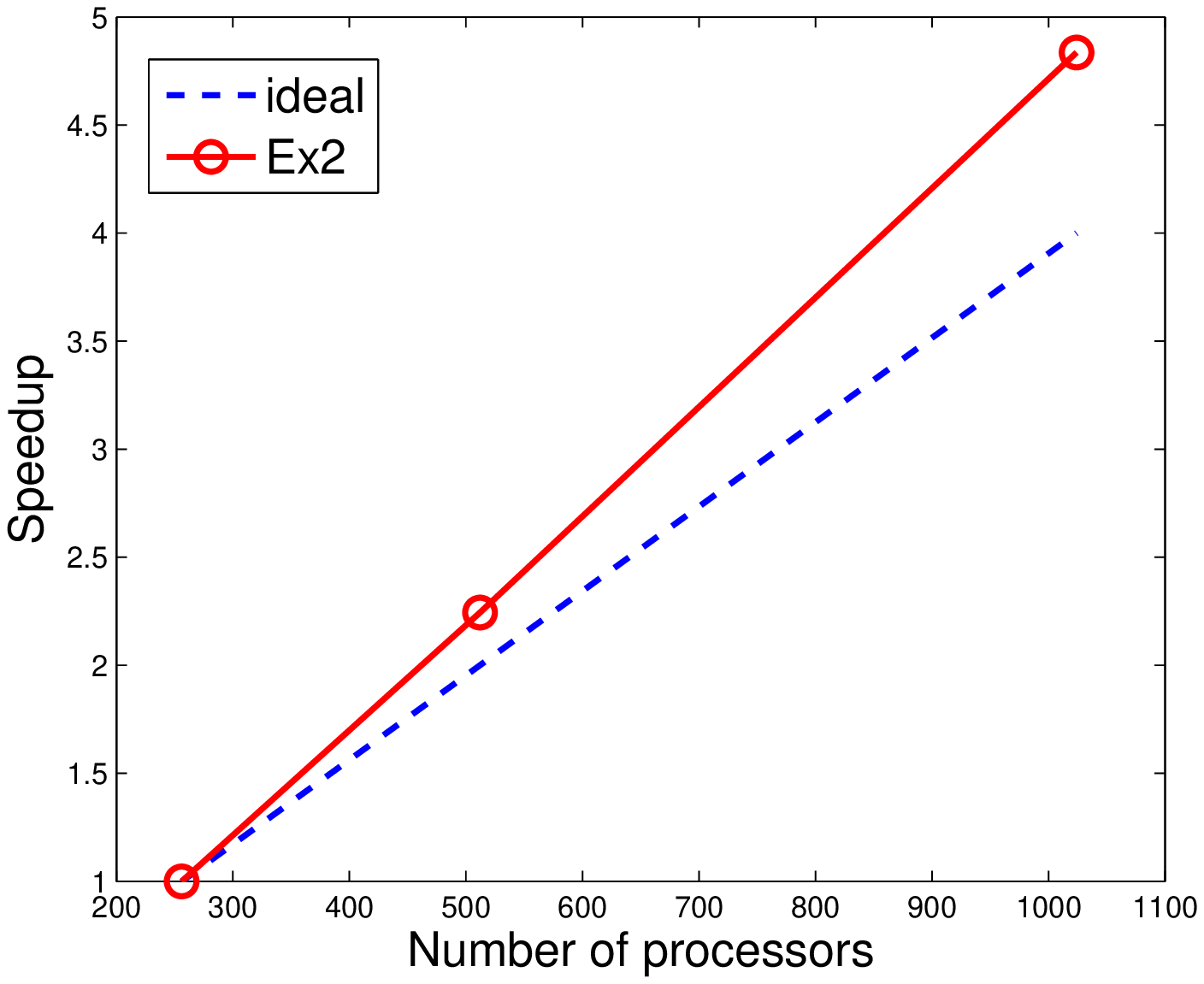}
\includegraphics[width=6cm]{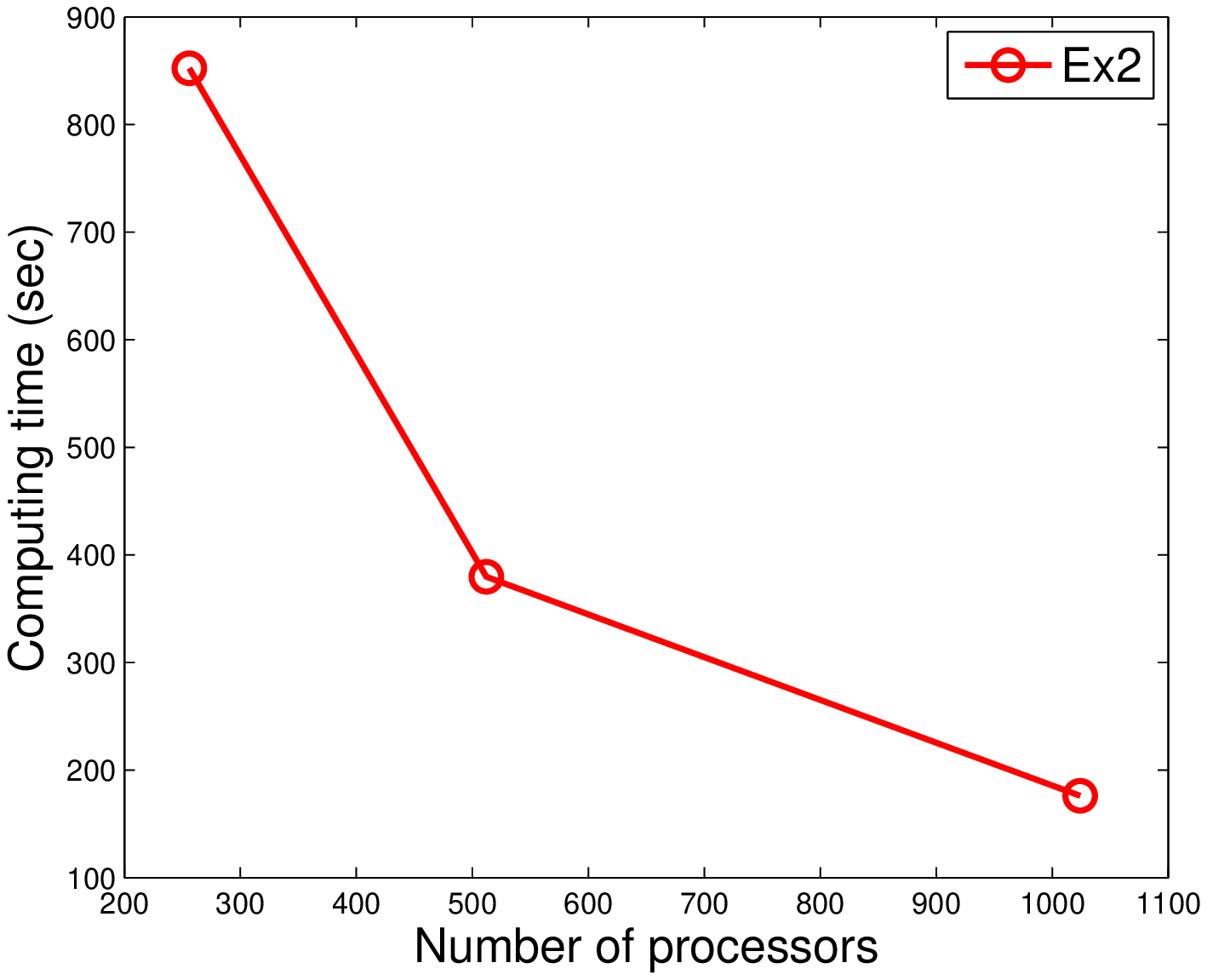}
\vspace{-0.15in}
\caption{The speedup (left) and computing time (right) for Example 2.}
\label{fig:scalex2}
\end{figure}

Using the same parameter settings as in Example 1, we perform the strong scalability test and the results are given in Table \ref{tab:mesh2} and Figure \ref{fig:scalex2}.
Superlinear speedup is obtained when $np\leq 1024$.
Next we test three sets of mesh and time step size in Table \ref{tab:ex2meshregu}. The $H^1$ regularization parameter is set to be $\beta_2=10^{-6}$, and 64 processors are used.
The overlap $iovlp=4$ and the fill-in level of ILU  $ilulevel=3$.
We observe from Table \ref{tab:ex2meshregu} that as the mesh and the time step size become finer,
the number of GMRES iterations grows slightly, and the computing time increases with the problem size.

\begin{figure}
\centering
\includegraphics[width=6cm]{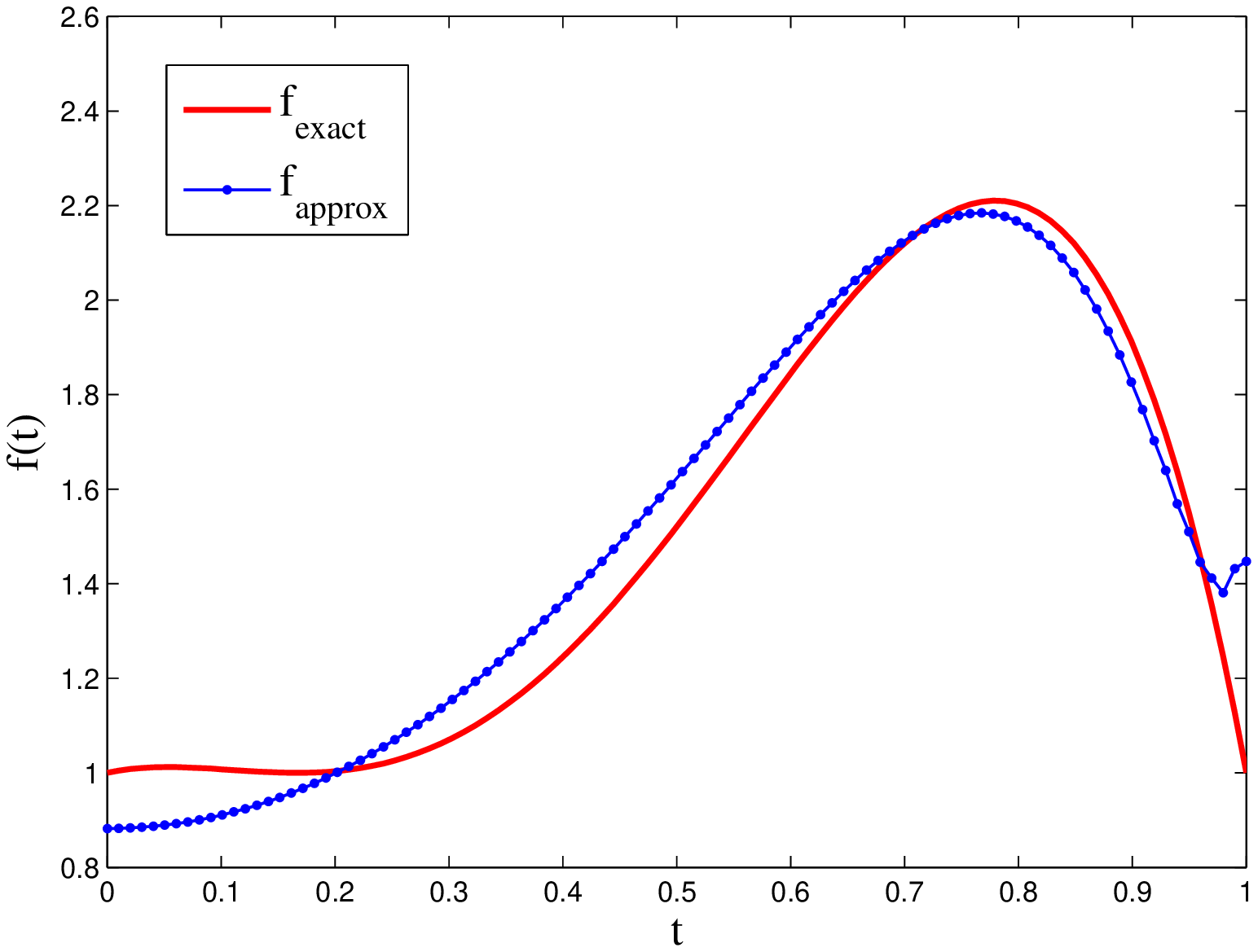}
\includegraphics[width=6cm]{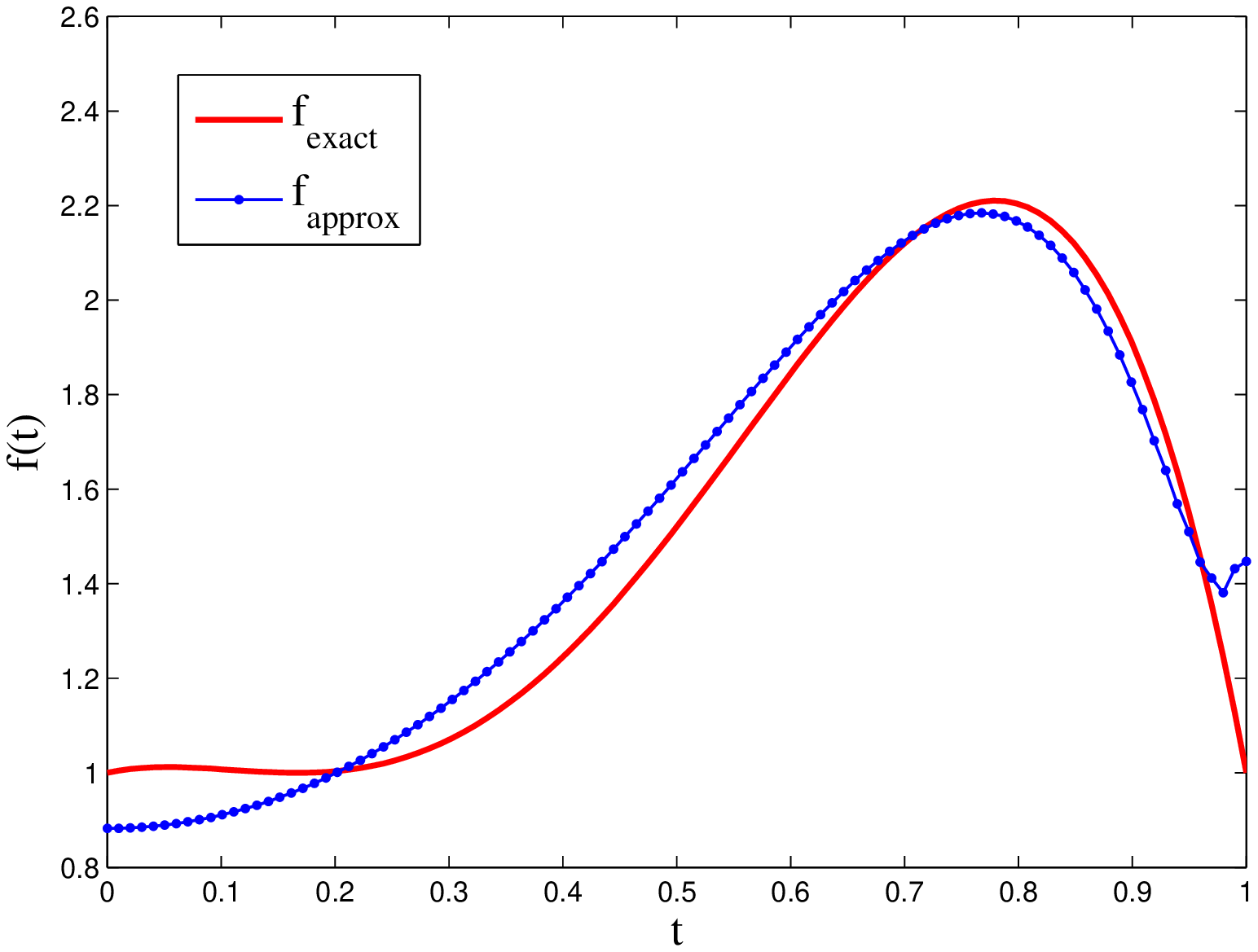}
\vspace{-0.15in}
\caption{Space-only ASM preconditioner (left) vs. space-time ASM preconditioner (right) for Example 2. }
\label{fig:pccomp}
\end{figure}

\begin{figure}
\centering
\includegraphics[width=6cm]{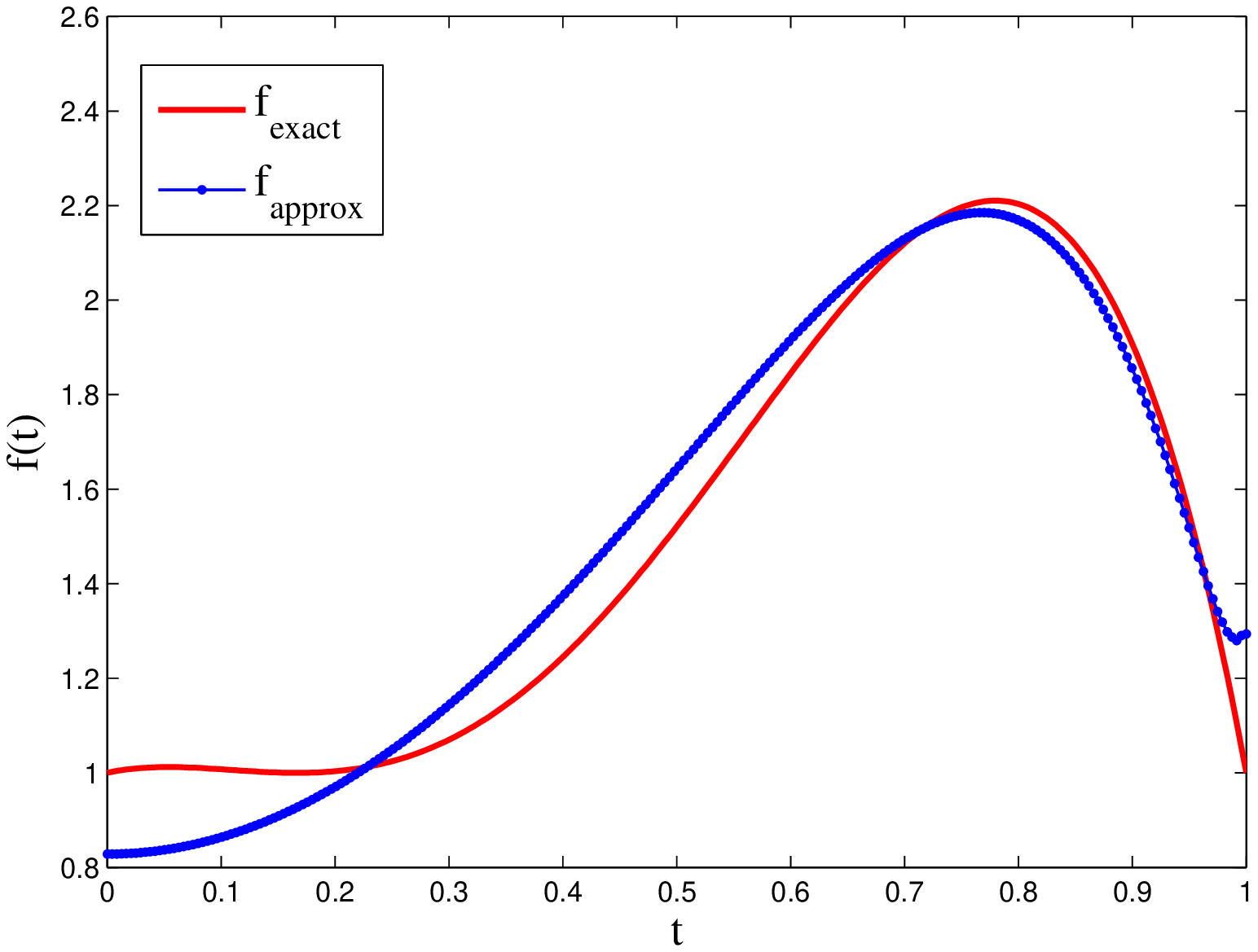}
\includegraphics[width=6cm]{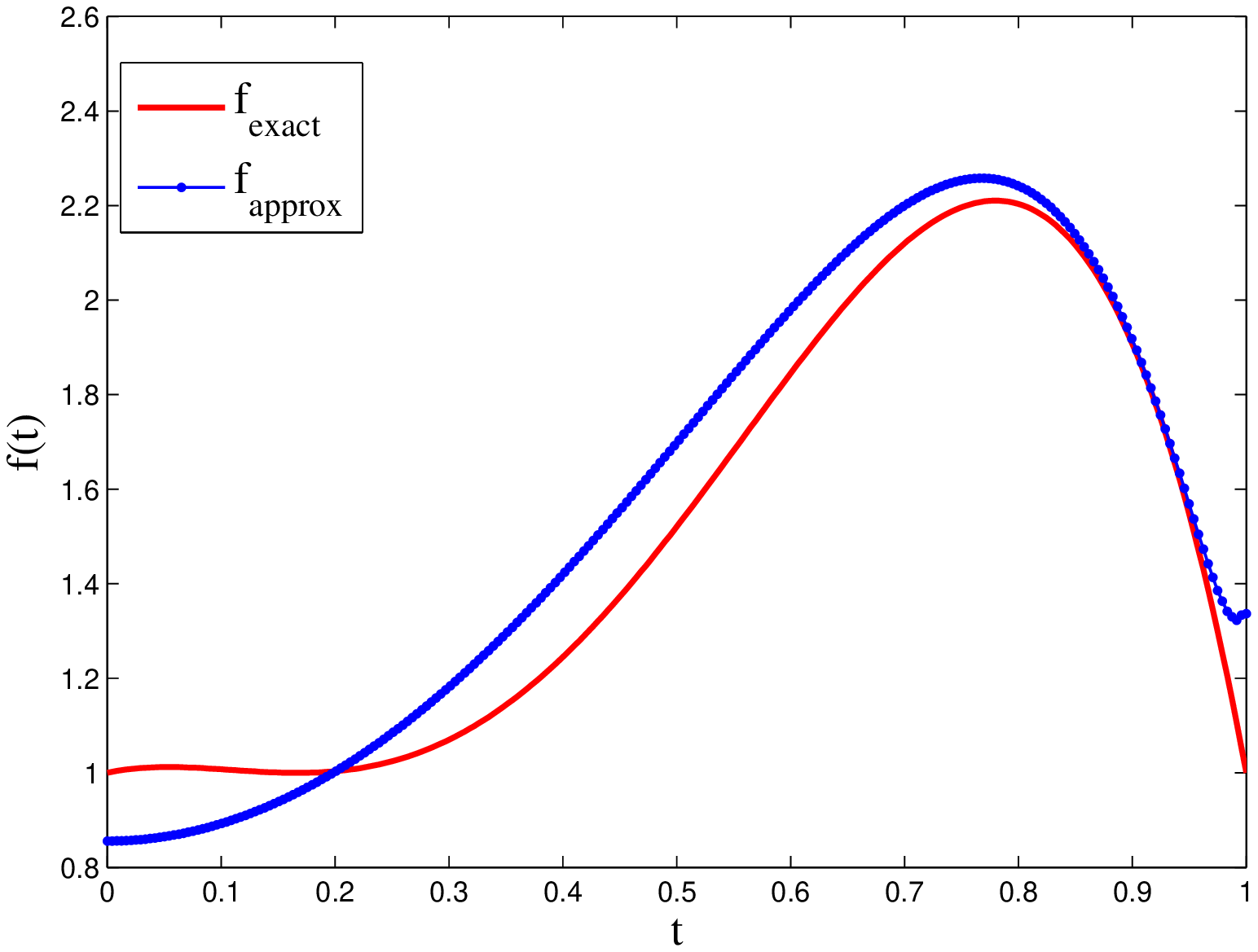}\\
\centering
\includegraphics[width=6cm]{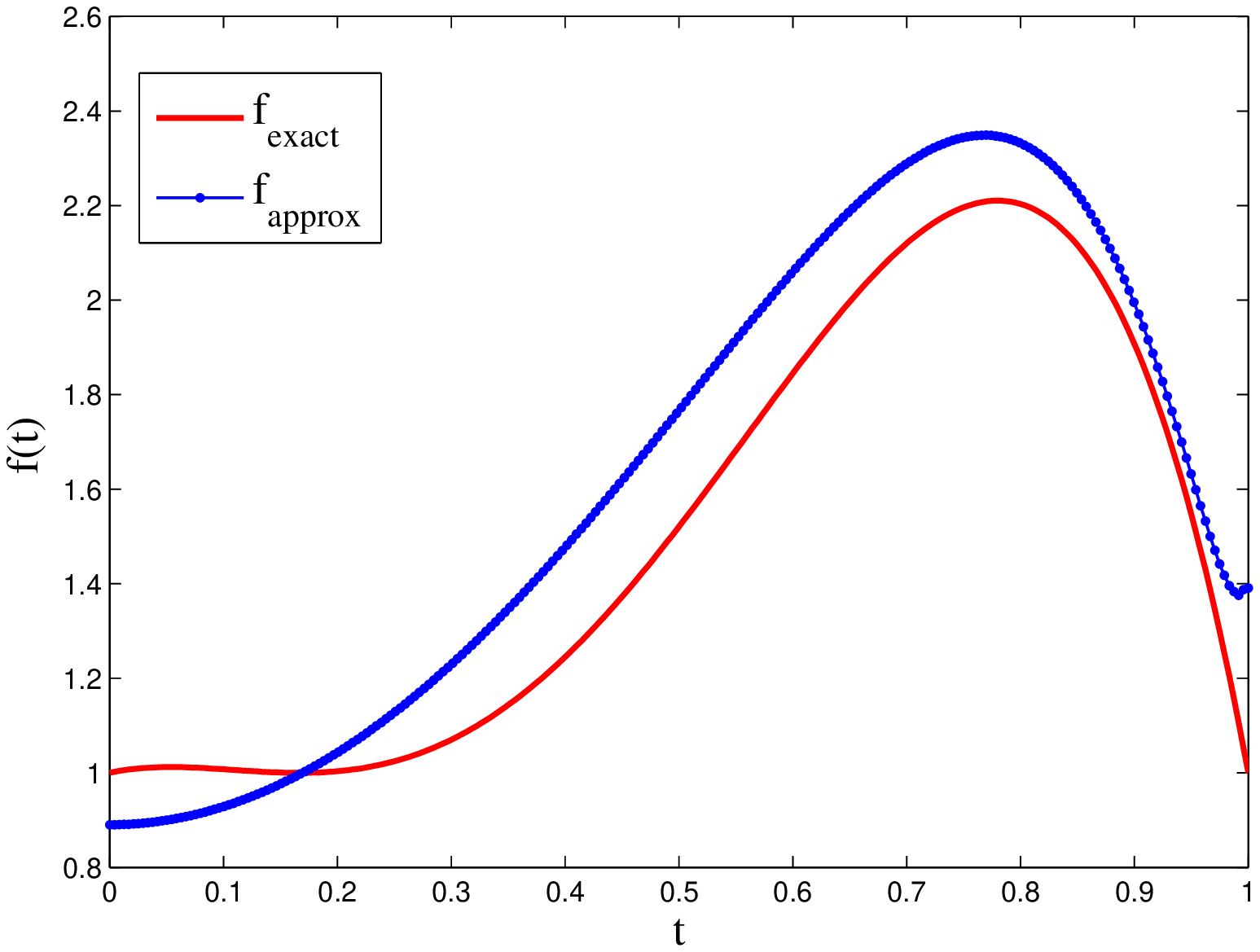}
\vspace{-0.15in}
\caption{Noise level test for example 2. Top left:  $\varepsilon=1\%$,  top right: $\varepsilon=5\%$,  bottom: $\varepsilon=10\%$. $n_t=240,~ n_x=64,~ n_y=64,~ DOF=1,966,320,~ \beta_2=10^{-5},~ np=64$.}
\label{fig:ex2noisetest}
\end{figure}

Now we investigate the performance of the space-time Schwarz preconditioner. An important feature of the proposed space-time Schwarz preconditioner lies in the parallelization in the time dimension.
If the time range is not partitioned as mentioned in the end of Section \ref{sec:space-time schwarz},
the preconditioner also works from the result in Figure \ref{fig:pccomp},
but it is observed from Table \ref{tab:pccomp}, under the same settings, that the ``space-only'' Schwarz preconditioner costs more iterations and computing time compared to the space-time Schwarz preconditioner.
Thus the Schwarz preconditioner with a partition in the time is more efficient than the ``space-only'' domain decomposition preconditioner.
In the end of this example, we perform the noise level test and the results are given in Figure \ref{fig:ex2noisetest}.
The results agree with our expectation that the reconstruction accuracy deteriorates with the increasing level of noise in the measurement data.
\begin{table}
\centering
\begin{tabular}{|ccccc|}
\hline
$np$&Its&Time(sec)&Speedup&Ideal\\
\hline
256&178&852.23&1&1\\
512&184&379.75&2.24&2\\
1024&247&176.23&4.84&4\\
\hline
\end{tabular}
\vspace{0.1in}
\caption{Scalability test for Example 2:  $n_t=160,~ n_x=160,~ n_y=160,~ DOF=8,192,160$.}
\label{tab:mesh2}
\end{table}
\begin{table}
\centering
\begin{tabular}{|ccccc|}
\hline
$n_t$&$n_x\times n_y$&DOF&Its&Time(sec)\\
\hline
100&$40\times 40$&320~100&40&8.44\\
240&$64\times 64$&1~966~320&46&33.25\\
320&$80\times 80$&4~096~320&49&69.28\\
\hline
\end{tabular}
\vspace{0.1in}
\caption{Mesh size and time step size test for Example 2: $\beta_2=10^{-6},~np=64$.}
\label{tab:ex2meshregu}
\end{table}
\begin{table}
\centering
\begin{tabular}{|c|ccc||ccc|}
\hline
Preconditioner&$np$&Its&Time(sec)&$np$&Its&Time(sec)\\
\hline
space-only&64&49&114.34&256&129&230.29\\
space-time&64&39&37.60&256&83&156.64\\
\hline
\end{tabular}
\vspace{0.1in}
\caption{Preconditioner comparison for Example 2:
$\beta_2=10^{-5}, ~n_t=100, ~n_x=40,~n_y=40,~DOF=320,100$ for $np=64$; $\beta_2=10^{-6}, ~n_t=320, ~n_x=80,~n_y=80,~DOF=4,096,320$ for $np=256$.}
\label{tab:pccomp}
\end{table}

\begin{figure}
  \centering
\includegraphics[width=6cm]{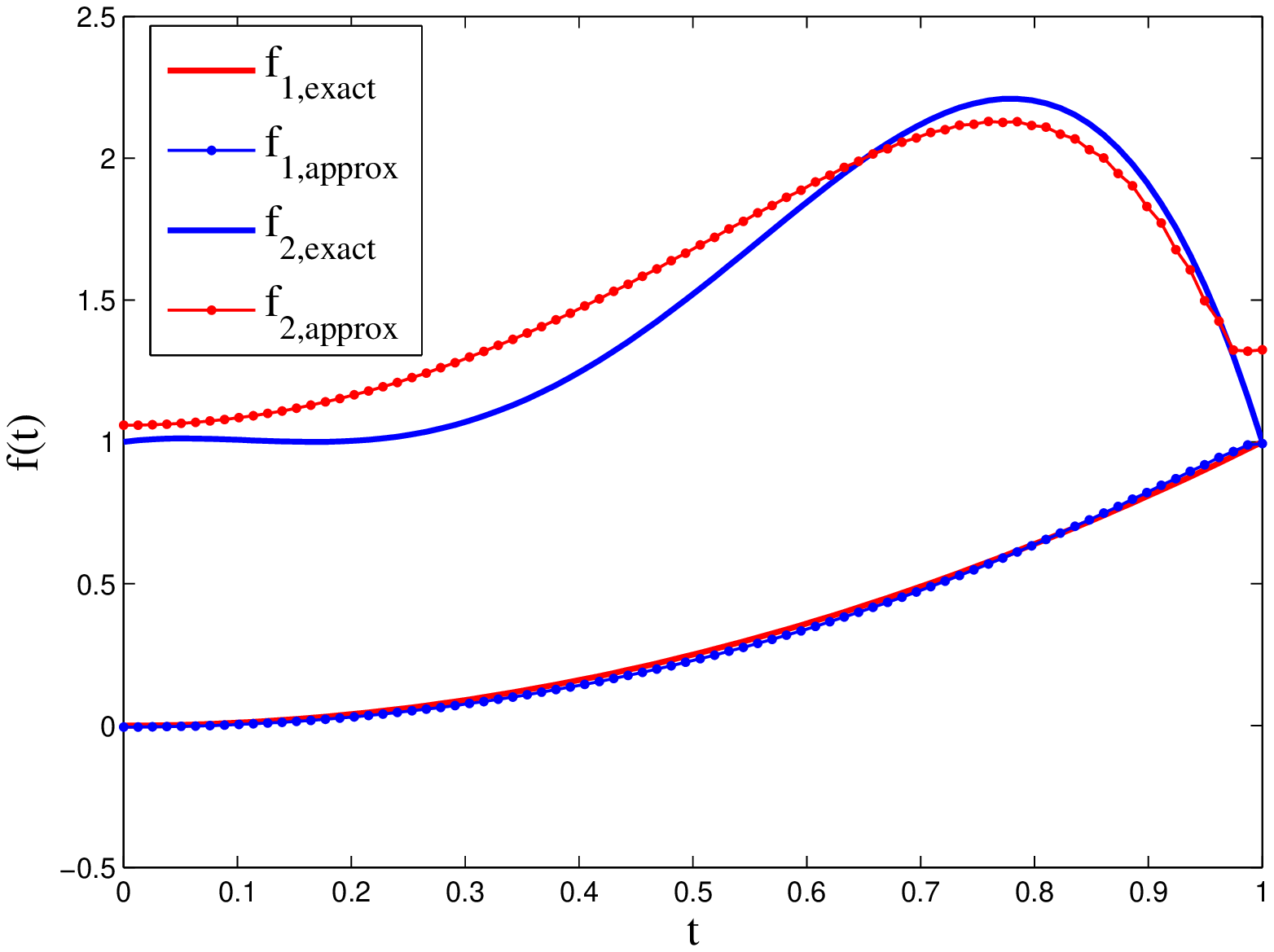}
\includegraphics[width=6cm]{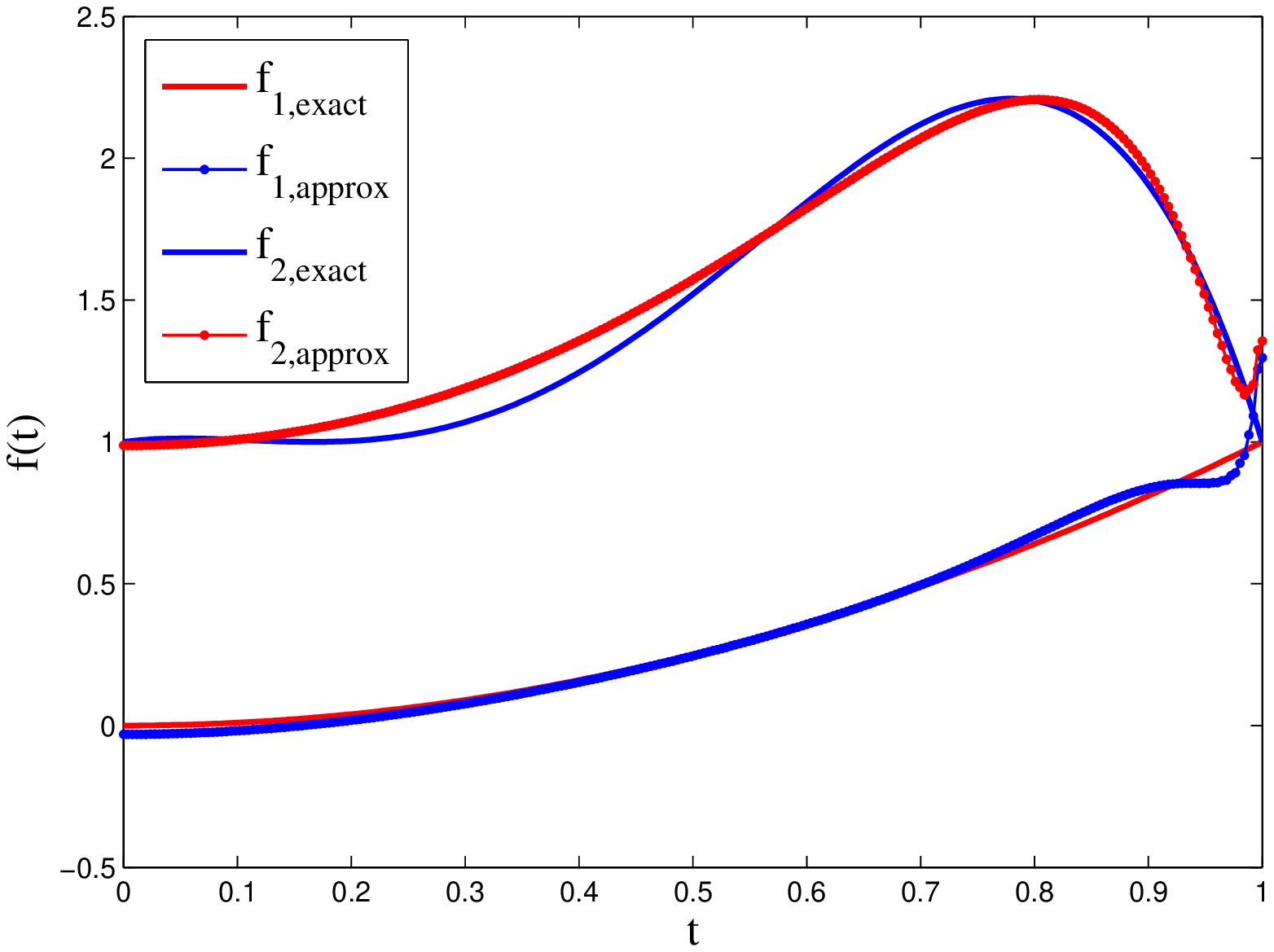}
\vspace{-0.15in}
\caption{Recovery of two source intensities using two sets of regularization parameters. Left: $\beta_1^1=0, \beta_1^2=10^{-4}, \beta_2^1=10^{-6}, \beta_2^2=10^{-5}$, right: $\beta_1^1=0, \beta_1^2=10^{-6}, \beta_2^1=0, \beta_2^2=10^{-6}$.}
\label{fig:2sources}
\end{figure}

\textbf{Example 3}. In this example, we test the recovery of two source intensities.
Compared with the single source intensity cases, more regularization parameters are needed.
Here we test the following two sets of regularization parameters with $np=64$:
\begin{enumerate}
\item[(1)]
$\beta_1^1=0, \beta_1^2=10^{-4}$ for the source intensity function $f_1$ and $\beta_2^1=10^{-6}, \beta_2^2=10^{-5}$ for $f_2$.
The mesh is $n_x=80, n_y=80$ and the time step is $n_t=80$;
\item[(2)]  $\beta_1^1=0, \beta_1^2=10^{-6}$ for $f_1$ and $\beta_2^1=0, \beta_2^2=10^{-6}$ for $f_2$.
The mesh and the time step are set to be $n_x=64, n_y=64$ and $n_t=256$.
\end{enumerate}
 The numerical results are shown in Figure \ref{fig:2sources}.
The computed $f_1$ matches with its original data perfectly, but the computed $f_2$ is less accurate.
As we see in the tests for Example 1 and Example 2, $f_2$ is physically harder to recover than the simpler function $f_1$.
Overall the reconstruction effect for both source intensities are reasonable.

\begin{figure}
\centering
\includegraphics[width=6cm]{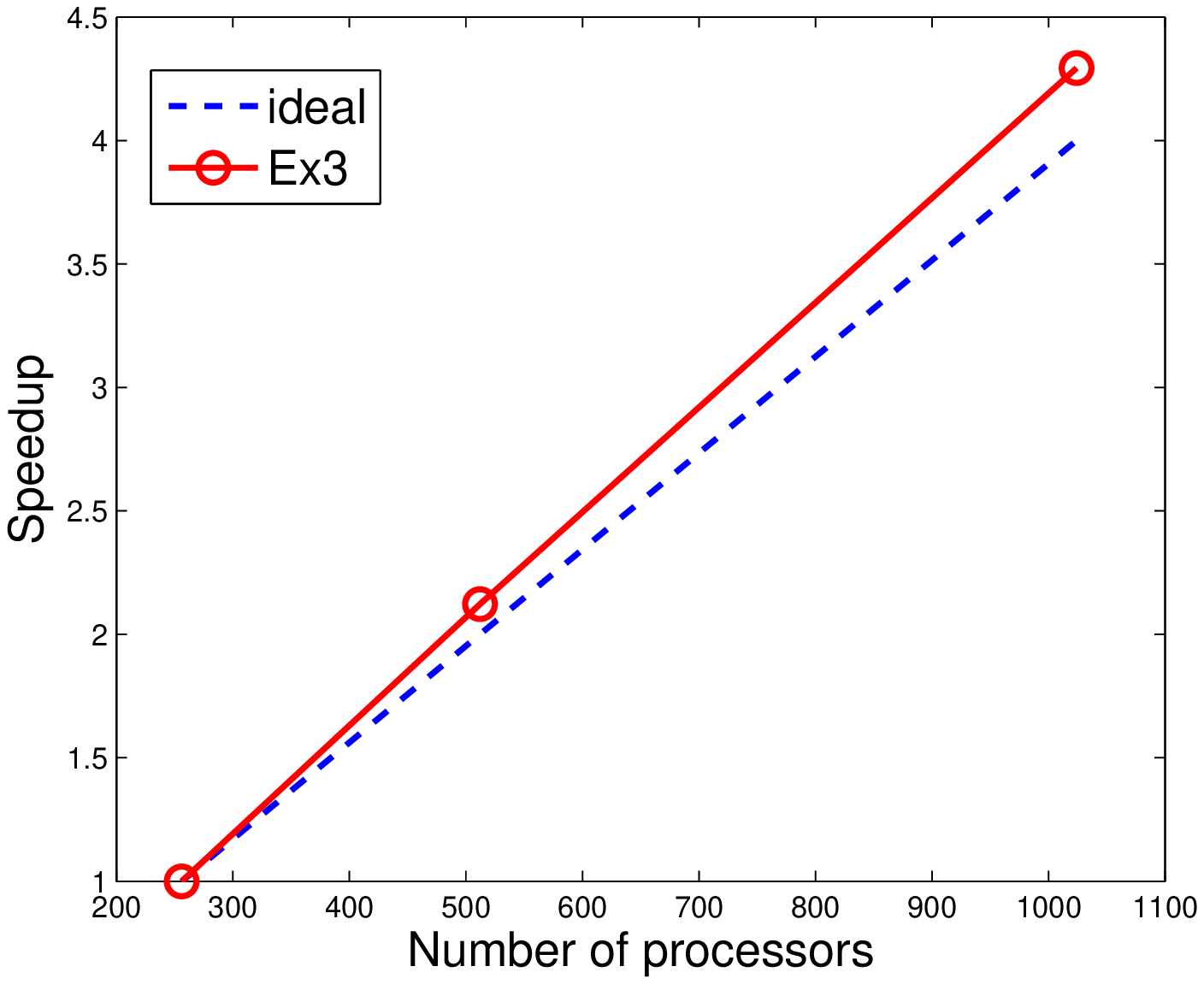}
\includegraphics[width=6cm]{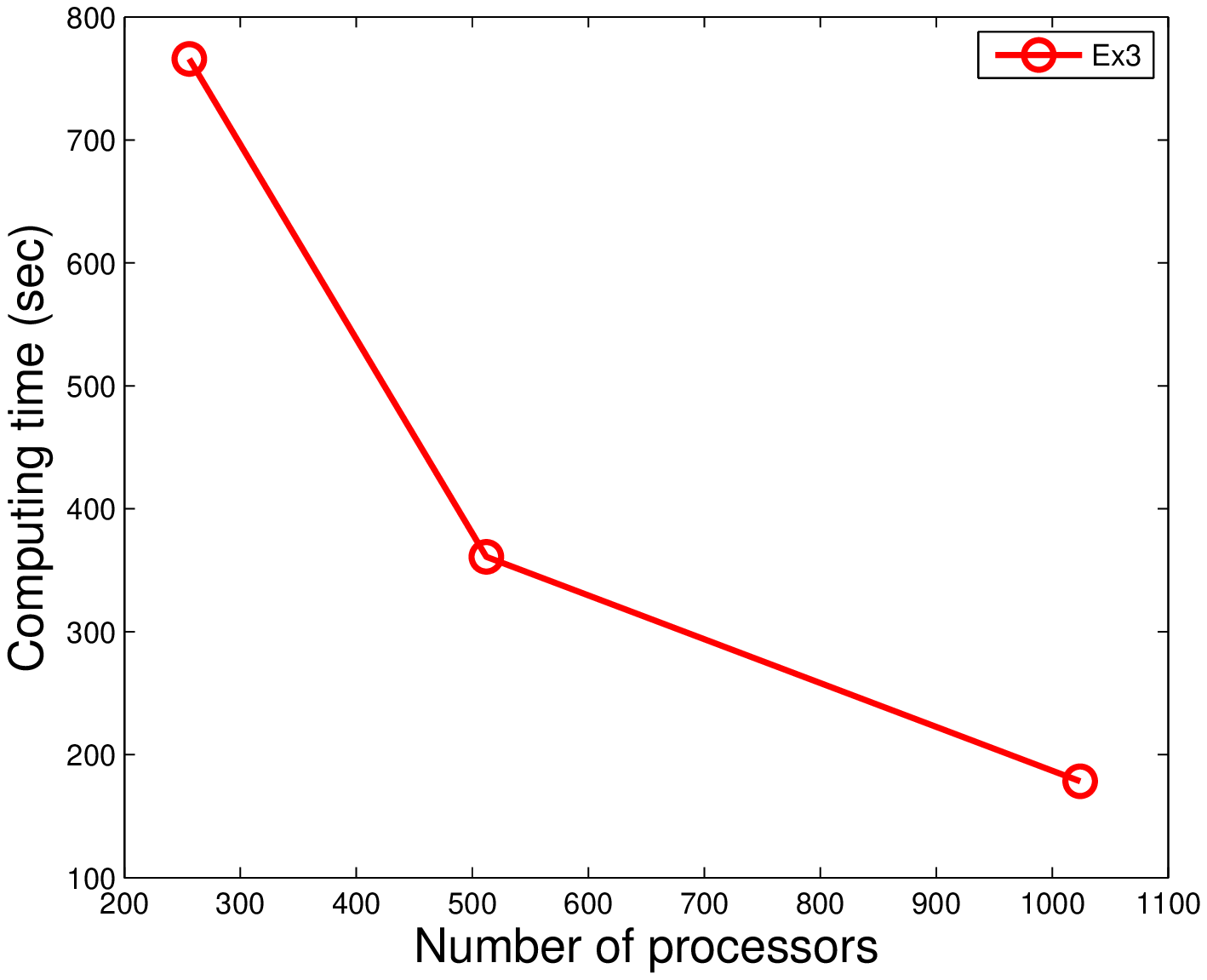}
\vspace{-0.15in}
\caption{The speedup (left) and the computing time (right) for Example 3.}
\label{fig:2sourcesscal}
\end{figure}

Next we show the strong scalability results in Figure \ref{fig:2sourcesscal} and Table \ref{tab:2sourcesscal}.
We still observe a superlinear speedup, although it is a bit worse than that of Examples 1 and 2.
It implies that it is more difficult to separate and identify multiple source intensities than the single source case.
And from our previous experiments with reduced space  SQP methods, the recovery of multiple sources is much more difficult to converge than that of the single source case.

Lastly, we test the algorithm with parameters such as the fill-in level of ILU factorization and the size of overlap in Table \ref{tab:2sourcesilulevel} and Table \ref{tab:2sourcesovlp}, respectively.

It is observed that the number of iterations decreases with the increase of the overlapping size or the fill-in level,
however, it costs more communication time when we increase the overlap between ``space-time'' subdomains,
and more computing time is used in the preconditioning stage when we raise the fill-in level of the ILU factorization.
\begin{table}
\centering
\begin{tabular}{|ccccc|}
\hline
$np$&Its&Time(sec)&Speedup&Ideal\\
\hline
256&148&765.80&1&1\\
512&150&360.88&2.12&2\\
1024&211&178.37&4.29&4\\
\hline
\end{tabular}
\vspace{0.1in}
\caption{Scalability test for  Example 3 with two point sources: $n_t=160,~ n_x=160,~ n_y=160,~ DOF=8,192,320$.}
\label{tab:2sourcesscal}
\end{table}

\begin{table}
\centering
\begin{tabular}{|ccc|}
\hline
$ilulevel$ &Its&Time(sec)\\
\hline
1&496&75.98\\
2&295&85.12\\
3&247&116.11\\
\hline
\end{tabular}
\vspace{0.1in}
\caption{Fill-in level of ILU test for Example 3: $\beta_1^2=10^{-6},~ \beta_2^2=10^{-6}$, $n_t=400,~ n_x=80,~ n_y=80,~ DOF=5,120,400,~np=128$.}
\label{tab:2sourcesilulevel}
\end{table}

\begin{table}
\centering
\begin{tabular}{|ccc|}
\hline
$iovlp$&Its&Time(sec)\\
\hline
1&-& -\\
2&432&114.22\\
4&247&121.26\\
6&238&400.12\\
\hline
\end{tabular}
\vspace{0.1in}
\caption{Overlap test for  Example 3:  $\beta_1^2=10^{-6},~ \beta_2^2=10^{-6}$, $n_t=400,~ n_x=80,~ n_y=80,~ DOF = 5,120,400,~np=128$.}
\label{tab:2sourcesovlp}
\end{table}

\textbf{Example 4}. We now test the numerical reconstruction for
three point sources, and observe how the speedup changes with increasing number of sources.

\begin{figure}
\centering
\includegraphics[width=8cm]{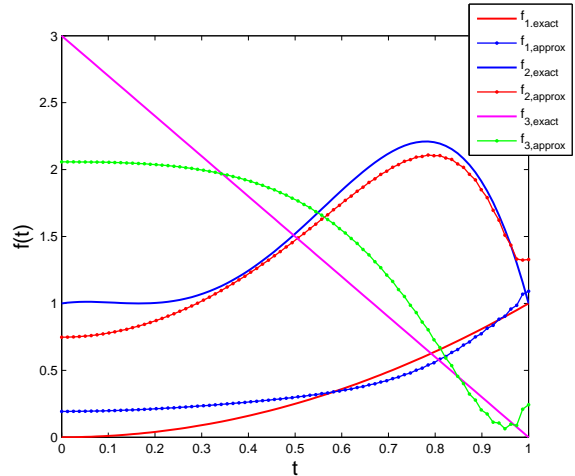}
\vspace{-0.15in}
\caption{The reconstruction for three point sources.}
\label{fig:ex4}
\end{figure}

\begin{figure}
\centering
\includegraphics[width=6cm]{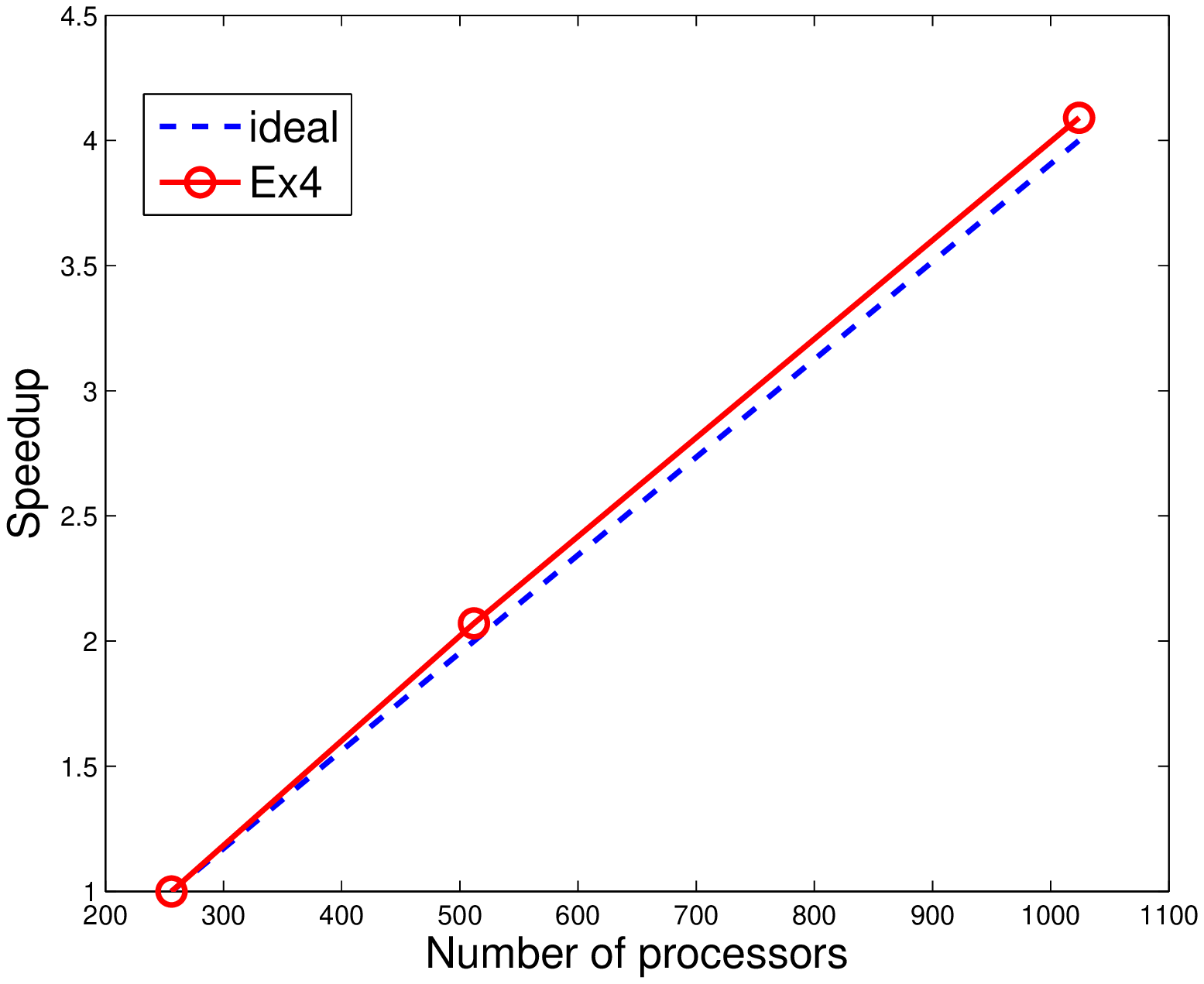}
\includegraphics[width=6cm]{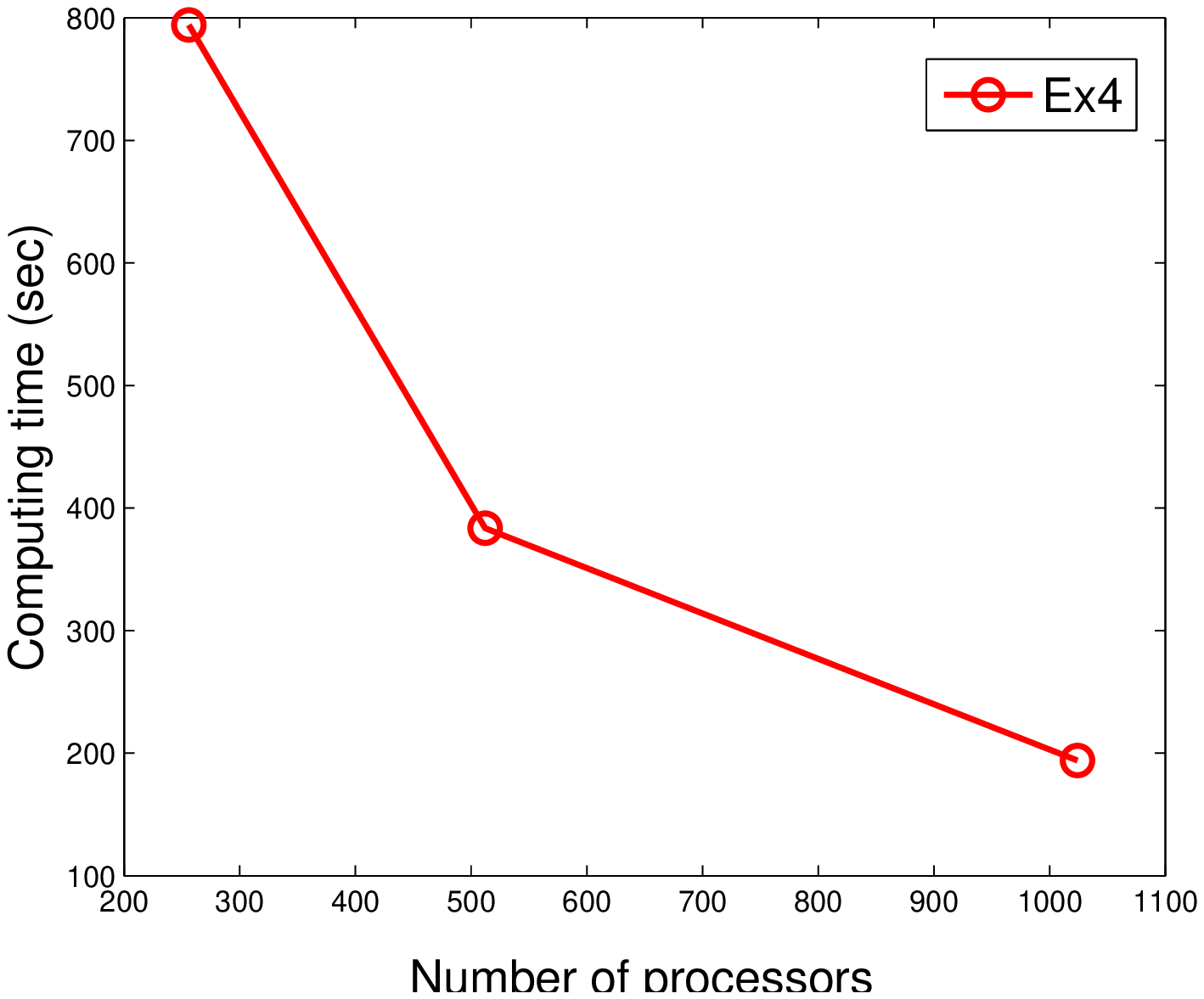}
\vspace{-0.15in}
\caption{The speedup (left) and the computing time (right) for Example 4.}
\label{fig:ex4scal}
\end{figure}

For this test, we take the spatial mesh $160\times 160$ and the number of time steps 160. Regularization parameters are respectively set to be $\beta_1^1=0, \beta_2^1=10^{-5}$, $\beta_1^2=10^{-4}, \beta_2^2=10^{-5}$ and $\beta_1^3=10^{-7}, \beta_2^3=8\times10^{-6}$.
From Figure \ref{fig:ex4} we see that, apart from the initial part of the third intensity which is not quite close to the true values, the rest are recovered satisfactorily.
Now we use a $160\times 160$ space mesh and $160$ time step to test the strong scalability and compute time
in Figure \ref{fig:ex4scal} and Table \ref{tab:3sourcesscal}. LU factorization is used as the subdomain solver.
It is observed that the speedup for three point sources is almost linear, still satisfactory but a bit worse than Examples 1,2 and 3.
As a conclusion the speedup deteriorates slowly with the number of unknown point sources.
\begin{table}
\centering
\begin{tabular}{|ccccc|}
\hline
$np$&Its&Time(sec)&Speedup&Ideal\\
\hline
256&148&794.19&1&1\\
512&150&383.56&2.07&2\\
1024&211&194.14&4.09&4\\
\hline
\end{tabular}
\vspace{0.1in}
\caption{Scalability test for  Example 4 with three point sources: $n_t=160,~ n_x=160,~ n_y=160,~ DOF=8,192,480$.}
\label{tab:3sourcesscal}
\end{table}
\subsection{Comparisons with two reduced space SQP methods}
A reduced space method for reconstruction of the location and intensity of a single point pollutant source was
developed in \cite{DZZ12}. With the source location known in our current case, the process of reconstructing the source intensity described in \cite{DZZ12} can be stated as follows:

\smallskip
\begin{algorithm}[H]
\caption{ ~~Nonlinear CG method}
\label{alg:nonlinearcg}
\begin{algorithmic}
      \STATE     Select the initial guesses $f^0$, and set $k:=0$.
      \STATE     For $k=1, 2, \cdots, N_{max}$
      \STATE      $\quad$ Solve the state system (the first equation in (\ref{eq:kkt2})) for $\{C_h^{n}(f^k)\}$;
      \STATE      $\quad$ Solve the adjoint system (the second equation in (\ref{eq:kkt2})) for $\{G_h^{n}(f^k)\}$;
      \STATE      $\quad$ Apply the nonlinear CG method to update $f^k$: $f^{k+1}=f^{k}+\alpha_1^k d^k$;
       \STATE      $\quad$ Stop the iteration if the stopping criteria are satisfied; otherwise set $k:=k+1$.
\end{algorithmic}
\end{algorithm}

We use the Fletcher-Reeves (FR) formula to update
the nonlinear CG direction $d^k$:
$d^k=J'_k+\gamma_k d^{k-1}$, with $d^0=J'_0$ and $\gamma_k={\|J'_k\|^2}/{\|J'_{k-1}\|^2}$ and
$J'_k=-(J_h^{\tau})'(f^k)$ being the negative gradient direction
which is obtained from formula (\ref{eq:Jwithf}). That is,
\[J'_k=-\left(\int_0^T (\beta_1 f^k(t) -  G_h^{\tau}(\x^*,t))g^{\tau}(t) d t +  \beta_2 \int_0^T  (f^k)'(t)(g^{\tau})'(t) d t\right)\]
We select the stepsize $\alpha_1^k$ such that
$\alpha_1^k= \mbox{argmin}_{\gamma>0} J_h^\tau(f^{k}+\gamma d_k)$.
For the $L^2$ and $H^1$ regularizations in (\ref{eq:reguterm}), we can work out the
exact formulae:
\beqnx
 \alpha_1^k&=&-\frac{(C_h^{M}(f^{k})-C^{\epsilon},A_h^{M})+\beta_1 (f^{k},d^k)}{ (A_h^{M},A_h^{M})+\beta_1 (d^k,d^k)}~~\mbox{($L^2$ regularization)},\\
 \alpha_1^k&=&-\frac{(C_h^{M}(f^{k})-C^{\epsilon},A_h^{M})+\beta_2 ((f^{k})',(d^k)')}{ (A_h^{M},A_h^{M})+\beta_2 ((d^k)', (d^k)')}~~\mbox{($H^1$ regularization)},
\eqnx
 where $A_h^{M}=C_h^{M}(f^{k})'d^{k}$ is obtained by solving the following sensitivity equation,
 \begin{equation*}
(\partial_{\tau} A_h^n, v_h) + (a\nabla  \bar{A}_h^n, \nabla v_h) + (\nabla\cdot (\mathbf{v}\bar{A}_h^n),v_h) \\
= v_h(\x^*) (\bar{d^k})^n, ~~\forall\,v_h\in \mathring{V}^h,~ n=1,\cdots,M,\\
\label{eq:sensitivity}
\end{equation*}
with $A_h^0=\mathbf{0}$. At each iteration, three time-dependent subsystems are solved.
When we implement this nonlinear CG method on parallel computers, we need to develop a
parallel solver for each subsystem. We will test two cases:
the first one uses the space domain decomposition
preconditioner but keeps the time marching process, while the second one uses a space-time
domain decomposition preconditioner as it is developed for the fully coupled system in this work and
solves each time-dependent subsystem all-at-once.
These two parallel solvers are denoted by RS(1) and RS(2) respectively.
We shall compare the computing times between our proposed space-time preconditioning method,
denoted by FS, and the two reduced space SQP methods RS(1) and RS(2); see
Table \ref{tab:compare}.
We use the three aforementioned methods to implement Example 2 with four sets of
meshes, and the number of processors increases with the refinement of the
meshes.
The subdomain solvers for all three kinds of methods are ILU.
We firstly compute the result by the FS method with zero initial guess and record the error accuracy
$e=\|f-f^*\|$. Then we use the same initial guess and the error bound $e$ for the reduced space methods RS(1) and RS(2),
and set the stopping criterium as $\|f^k-f^*\|<e$. In this way we can compare the computing time for all these methods.

As shown in Table \ref{tab:compare}, the computing time of the FS method is
much less than the ones of RS(1) and RS(2).
For the  two reduced space methods, RS(2) using the space-time domain decomposition
solver is faster than RS(1) keeping the time marching process and using a space domain decomposition solver.
We can see that the space-time fully coupled preconditioner is much better for parallellization,  and
the all-at-once method for the fully coupled KKT system is always more efficient than the reduced space
iterative optimization method on parallel systems.

\begin{table}
\centering
\begin{tabular}{|ccccc|}
\hline
$np$&$n_t$&$n_x\times n_y$&Solver&Time(sec)\\
\hline
64&40&$40\times40$&FS&12.064\\
&&&RS(1)&418.580\\
&&&RS(2)&125.484\\
\hline
128&80&$80\times80$&FS&15.525\\
&&&RS(1)&682.794\\
&&&RS(2)&99.528\\
\hline
256&160&$80\times80$&FS&23.736\\
&&&RS(1)&994.962\\
&&&RS(2)&200.543\\
\hline
512&320&$160\times160$&FS&136.717\\
&&&RS(1)&7240.881\\
&&&RS(2)&1094.886\\
\hline
\end{tabular}
\vspace{0.1in}
\caption{The computing time of the proposed full space method FS, the reduced space method RS(1) and RS(2).}
\label{tab:compare}
\end{table}
\begin{section}{Concluding remarks}\label{sec:concluding remarks}
We developed a new space-time domain decomposition method for unsteady source inversion problems.
The main ingredient of our algorithm includes solving the fully coupled KKT system by GMRES iteration
with a space-time additive Schwarz preconditioner.
Although the size of the linear system is significantly increased compared to the reduced space SQP methods, the one-shot method avoids the sequential step between the state equation and the adjoint equation, as well as the time-marching process in the time dimension, and thus achieves higher degree of parallelism.
This is well confirmed by the numerical results shown in the last section.
Another advantage of the new method is that the recovery of multiple sources is obtained using the same algorithmic and software framework as the single source case,
and the framework is easily extended to recover other kinds of source intensities.

We have observed from the numerical examples that the new space-time additive Schwarz method
is quite robust also with respect to  the noise in the observation data.
It is important to note that the new space-time method is highly parallel and scalable, and
extensible naturally to three-dimensional problems.
\end{section}

\section*{Acknowledgements}
The authors would like to thank two anonymous referees for their many insightful and
constructive comments and suggestions, which have led to a great improvement of the results and organisation
of this work.

\end{document}